\documentclass{gtart_sp}   
\pdfoutput=1

\usepackage[matrix,arrow,curve,frame]{xy}       

\title{On the chain-level intersection pairing for PL manifolds}

\author{J\,E\,McClure}
\givenname{J\,E\,}
\surname{McClure}
\address{Department of Mathematics\\
Purdue University\\\newline
150 N\,University Street\\
West Lafayette, IN  47907-2067\\USA}
\email{mcclure@math.purdue.edu}
\urladdr{}

\volumenumber{10}
\issuenumber{}
\publicationyear{2006}
\papernumber{35}
\lognumber{0508}
\startpage{1391}
\endpage{1424}

\doi{}
\MR{}
\Zbl{}

\arxivreference{math.QA/0410450}


\received{22 Oct 2004}
\revised{1 August 2006}
\accepted{14 July 2005}
\published{4 October 2006}
\publishedonline{4 October 2006}
\proposed{Gunnar Carlsson}
\seconded{Haynes Miller, Bill Dwyer}
\corresponding{}
\editor{Ralph Cohen}
\version{}

\let\xysavmatrix\xymatrix
\def\xymatrix{\disablesubscriptcorrection\xysavmatrix}
\AtBeginDocument{\let\bar\wbar\let\tilde\wtilde\let\hat\what}
\makeop{Ch}
\newcommand{\ov}{\wbar}

\xymatrixcolsep{1.9pc}                          
\xymatrixrowsep{1.9pc}
\newdir{ >}{{}*!/-5pt/\dir{>}}                  

\makeatletter
\def\cnewtheorem#1[#2]#3{\newtheorem{#1}{#3}[section]
\expandafter\let\csname c@#1\endcsname\c@theorem}

\newtheorem{theorem}{Theorem}[section]
\cnewtheorem{proposition}[theorem]{Proposition}
\cnewtheorem{lemma}[theorem]{Lemma}
\cnewtheorem{conjecture}[theorem]{Conjecture}
\cnewtheorem{corollary}[theorem]{Corollary}

\theoremstyle{definition}
\cnewtheorem{definition}[theorem]{Definition}
\cnewtheorem{pdef}[theorem]{Provisional Definition}
\cnewtheorem{example}[theorem]{Example}
\cnewtheorem{notation}[theorem]{Notation}
\cnewtheorem{question}[theorem]{Question}
\cnewtheorem{construction}[theorem]{Construction}
\cnewtheorem{assumption}[theorem]{Assumption}

\theoremstyle{remark}
\cnewtheorem{remark}[theorem]{Remark}
\makeatother  
\makeautorefname{pdf}{Provisional Definition}

\DeclareMathOperator{\Map}{Map}

\DeclareMathOperator*{\colim}{colim}

\newcommand{\F}{{\cal F}}

\renewcommand{\a}{{\mathbf a}}

\renewcommand{\O}{{\cal O}}

\newcommand{\supp}{{\text{\rm supp}}}
\renewcommand{\int}{{\text{\rm int}}}
\newcommand{\im}{{\text{\rm im}}}

\newlength{\labwidth}
\newcommand{\labarrow}[1]{
\settowidth{\labwidth}{$\scriptstyle \;\; #1 \;\;$}
\stackrel{#1}{\smash{\hbox to \labwidth{\rightarrowfill}} 
\vphantom{\longrightarrow}}
}

\begin{document}

\begin{asciiabstract}
Let M be a compact oriented PL manifold and let C_*M be its PL
chain complex.  The domain of the chain-level intersection pairing is a
subcomplex G of C_*M\otimes C_*M.  We prove that G is a ``full''
subcomplex, that is, the inclusion of G in C_*M \otimes C_*M is a
quasi-isomorphism.  An analogous result is true for the domain of the
iterated intersection pairing.  Using this, we show that the intersection
pairing gives C_*M a structure of partially defined commutative DGA, which
in particular implies that C_*M is canonically quasi-isomorphic to an
E_\infty chain algebra.
\end{asciiabstract}

\begin{htmlabstract}
Let M be a compact oriented PL manifold and let C<sub>*</sub>M be its PL
chain complex.  The domain of the chain-level intersection pairing is a
subcomplex of C<sub>*</sub>M&otimes; C<sub>*</sub>M.  We prove that
the inclusion map from this subcomplex to C<sub>*</sub>M &otimes; C<sub>*</sub>M is a
quasi-isomorphism.  An analogous result is true for the domain of the
iterated intersection pairing.  Using this, we show that the intersection
pairing gives C<sub>*</sub>M a structure of partially defined commutative DGA, which
in particular implies that C<sub>*</sub>M is canonically quasi-isomorphic to an
E<sub>&</sub>infin; chain algebra.
\end{htmlabstract}

\begin{abstract}
Let $M$ be a compact oriented PL manifold and let $C_*M$ be its PL
chain complex.  The domain of the chain-level intersection pairing is a
subcomplex of $C_*M\otimes C_*M$.  We prove that 
the inclusion map from this subcomplex to $C_*M \otimes C_*M$ is a
quasi-isomorphism.  An analogous result is true for the domain of the
iterated intersection pairing.  Using this, we show that the intersection
pairing gives $C_*M$ a structure of partially defined commutative DGA, which
in particular implies that $C_*M$ is canonically quasi-isomorphic to an
$E_\infty$ chain algebra.
\end{abstract}
\maketitle
{\leftskip 25pt\small\hyperlink{Err}{Erratum attached}\par}

\section{Introduction}

\label{sec1}

Let $M$ be a compact oriented PL manifold. The chain-level intersection 
pairing was introduced by Lefschetz \cite{Lef} as a tool for constructing 
the intersection pairing on the homology of $M$. 
A version of the chain-level intersection pairing is a 
basic ingredient in Chas and Sullivan's construction \cite{CS} of a 
Batalin--Vilkovisky structure on the homology of the free loop space of $M$.

For a complete understanding of the chain-level intersection pairing, it seems 
helpful to have the following theorem.  Let $C_*M$ be the PL chain complex 
of $M$ (see \fullref{sec2} for the definition). 
Let us say that a subcomplex of a chain complex is {\it full\/} if the 
inclusion map is a quasi-isomorphism.  

\begin{theorem}
\label{1.1}
The domain of the chain-level intersection pairing is a full
subcomplex of $C_*M\otimes C_*M$.
\end{theorem}

It might seem at first that something like \fullref{1.1} would have been 
needed already by Lefschetz to define the intersection pairing on homology,
but for that purpose two weaker facts would suffice:

(i)\qua
For any cycles $C$ and $D$ in $C_* M$, 
the chain $C\otimes D$ is homologous to an element in the domain of the intersection
pairing.

(ii)\qua
If $C',D'$ are two other cycles with $C'\otimes D'$ homologous to $C\otimes
D$, then the difference $C'\otimes D'-C\otimes D$ is the boundary of an element in the
domain of the intersection pairing. (This is needed to show that the
intersection pairing on homology is well-defined.)

\fullref{1.1} is harder to prove than (i) and (ii) because (among other
reasons) a cycle in $C_*M\otimes C_* M$ cannot in general be written in the 
form $\sum C_i\otimes D_i$ with $C_i$ and $D_i$ cycles.

One goal of this paper is to prove \fullref{1.1} and, more
generally, the analogous statement for the $k$--fold iterate of 
the intersection pairing; see \fullref{11.b.3} and \fullref{11.b.4}.

It seems useful to go farther and to show that the intersection pairing gives 
$C_*M$ a structure of ``partially defined commutative DGA;''  this is 
the second (and main) goal of this paper (see \fullref{11.1}).  
Combining this with \fullref{rev18} and \cite[Theorem 1]{Wilson} 
shows in particular that $C_*M$ is canonically quasi-isomorphic to an 
$E_\infty$ chain algebra.

The third goal of this paper is to give a new treatment of the
chain-level intersection pairing, based on the account of Goresky and MacPherson \cite{GM} but
with some improvements. 

The results of this paper will be applied by the author \cite{McClure} to prove 
two theorems about the Chas--Sullivan operations.  Let $LM$ be the free loop
space of $M$, let $S_*$ denote the singular chain functor and let $\F$ be the 
framed little 2-disks operad as in \cite{Getzler}.

\medskip
{\bf Theorem A}\qua 
{\sl
The Batalin--Vilkovisky structure on the homology of $LM$ is induced by a
natural action of an operad quasi-isomorphic to $S_*\F$ on a chain complex
quasi-isomorphic to $S_*(LM)$.
}
\medskip

(Theorem A is the analog for $H_*(LM)$ of Deligne's Hochschild cohomology 
conjecture; see Markl, Shnider and Stasheff \cite[Section I.1.19]{MSS}.)

\medskip
{\bf Theorem B}\qua
{\sl
The Eilenberg--Moore spectral sequence converging to the homology of $LM$ is
a spectral sequence of Batalin--Vilkovisky algebras.
}
\medskip

The paper is organized as follows.

\fullref{sec1a} gives a brief discussion of the definitions of the 
chain-level intersection pairing given in \cite{Lef}, \cite{Lef2} and 
\cite{GM} and 
explains why these versions of the definition
are not convenient as a starting point for proving \fullref{1.1}.

\fullref{sec2} recalls (from \cite{GM}) the definition of the PL chain 
complex of a PL space.
\fullref{sec3} recalls (also from \cite{GM}) a method for making chain-level
constructions by using relative homology.
\fullref{sec4} constructs the umkehr (that is, ``reverse'') map in 
relative homology induced by a PL map between compact oriented PL manifolds.  
In \fullref{sec5} a chain-level umkehr map is deduced from this 
using the method of \fullref{sec3}.
\fullref{sec6} recalls the definition of exterior product for PL
chains.

In \fullref{sec7}, the chain-level intersection pairing is defined as 
the composite of the exterior product and the chain-level umkehr map 
induced by the diagonal map; the motivation for this definition is that 
the intersection of two subsets of a set $S$ can be identified with the 
intersection of their Cartesian product with the diagonal in $S\times S$.

\fullref{sec8} gives the formal definition of ``partially defined 
commutative algebra.'' I use Leinster's concept of homotopy algebra 
\cite{Leinster} for this purpose rather than the Kriz--May definition of 
partial algebra \cite{KM} (but I will use the term ``Leinster partial 
algebra'' instead of ``homotopy algebra,''  since the latter term seems 
excessively generic).  The reason for using Leinster's definition is that it is
simpler and more intuitive.  It will be shown in \cite{McClure} that the 
Kriz--May definition is a special case of the Leinster definition (see \fullref{rev2}(b) below).

\fullref{sec9}--\fullref{sec12} give the proof that the intersection
pairing and its iterates determine a Leinster partial commutative DGA 
structure on $C_*M$.  The proof uses a ``moving lemma'' (\fullref{12.b.1})
which is proved in \fullref{sec12a} by means of a general-position result (\fullref{rev5}) that may be of 
independent interest.
\fullref{rev5} is proved in \fullref{sec13} and 
\fullref{sec13a}.

My work on this paper was partially supported by NSF grant DMS-0405693.

I would like to thank Mike Mandell for his help, Greg Friedman, Mark Goresky 
and Clint McCrory for useful correspondence, and the referee for a careful 
reading of the paper and for helpful comments.  I would especially like to 
thank Shmuel Weinberger for referring me to \cite{GM}.

I would also like to thank the Lord for making my work possible.

\section{The Lefschetz and Goresky--MacPherson definitions of the
chain-level intersection pairing}

\label{sec1a}

This section is not needed logically for the rest of the paper; it is
offered as motivation for \fullref{sec2}--\fullref{sec7}.   The reader may
also find it helpful to consult Steenrod's account of Lefschetz's
work on the intersection pairing \cite[pages 28--30]{Steenrod}.

This section uses some technical terms which will be defined in
\fullref{sec2}--\fullref{sec6}.

Lefschetz's first account of the chain-level intersection pairing $C\cdot D$ 
was in \cite{Lef}.  In this paper he uses the obvious definition: 
if $C=\sum m_i\sigma_i$ and
$D=\sum n_i\tau_i$ then
\begin{equation}
\label{l}
C\cdot D
=\sum \pm\, m_in_j\,\sigma_i\cap \tau_j,
\end{equation}
where the signs are determined by the orientations of $\sigma_i$, $\tau_j$ and
$M$.  This formula does not in fact give a chain unless all of the 
intersections $\sigma_i\cap \tau_j$ have the same dimension, so some 
restriction on the pair $(C,D)$ is necessary.
Generically, the intersection of a $p$--dimensional PL subspace and a 
$q$--dimensional PL subspace
has dimension $\leq p+q-\dim M$; pairs of PL subspaces with this property 
are said to be in {\it general position\/}.  Lefschetz restricts the domain of 
the intersection pairing to pairs $(C,D)$ for which all of the pairs
$(\sigma_i,\tau_j)$ are in general position, and 
he interprets terms $\sigma_i\cap\tau_j$ which are in dimension less than 
$\dim C+\dim D-\dim M$ as 0.

In order to prove the crucial formula 
\begin{equation}
\label{1a.1}
\partial(C\cdot D)=(\partial C)\cdot D \,\pm\, C\cdot \partial D,
\end{equation}
Lefschetz  imposes a further restriction on the domain of the intersection 
pairing:  he requires that all of the pairs 
$(\partial\sigma_i,\tau_j)$ and $(\sigma_i,\partial\tau_j)$ 
should also be in general position.%
\footnote{If $C$ and $D$ are chains on the {\it same triangulation\/} this 
condition forces $C\cdot D$ to be 0, because
$\sigma_i$ and $\tau_j$ will intersect along a common face and therefore
$\sigma_i\cap\tau_j$ will be contained in $\partial\sigma_i\cap\tau_j$.}
This assumption allows him to prove 
equation \eqref{1a.1} by working with one pair of simplices at a time and 
extending additively.

This definition has the disadvantage that the domain of the
intersection pairing is not invariant under subdivision.  For example, if 
$\sigma$ and $\tau$ are $1$--simplices in a $2$--manifold intersecting at a 
point in the interior of both, then the pair $(\sigma,\tau)$ is in the domain, 
but if we subdivide $\sigma$ and $\tau$ at the intersection point we
obtain a pair of chains $(\sigma'{+}\sigma'',\tau'{+}\tau'')$ which is not in the
domain (because for example the pair $(\partial\sigma',\tau')$ is not in 
general position).%
\footnote{Note also that, if the intersection point is 
$P$, then formula \eqref{l} gives $\sigma\cdot\tau=\pm P$ but 
$(\sigma'+\sigma'')\cdot (\tau'+\tau'')=\pm 4P$.}
This phenomenon is general: if $(C,D)$ is in the domain 
of this version of the intersection pairing with $C\cdot D\neq 0$ then there 
will always be a subdivision in which the pair of chains determined by $C$ 
and $D$ is not in the domain. 

Lefschetz returned to the chain-level intersection pairing in \cite[Section
IV.6]{Lef2}.  He gave a formula more general than \eqref{l} (equation
(46) on page 212) in which the coefficients are ``looping coefficients''
\cite[Section IV.5]{Lef2}.
This allowed him to enlarge the domain of the intersection pairing as follows:
if we write $\supp(C)$ for the union of the simplices that occur in $C$, then
$C\cdot D$ is defined when the three pairs
$(\supp(C),\supp(D))$, $(\supp(\partial C),\supp(D))$, 
$(\supp(C),\supp(\partial D))$ are in general position; note that this 
condition is invariant under subdivision.

The ``looping coefficients'' used in Lefschetz's second 
definition are tricky to define explicitly \cite[Subsection 58 on page
216]{Lef2}.  The theory has been worked out carefully in Keller \cite{Keller} (which
I have not had an opportunity to consult) and seems to be rather complicated
(see the Math Review).

The chain-level intersection pairing became temporarily obsolete when the cup
product was discovered and it was noticed that the intersection pairing in 
homology could be defined using only Poincar\'e duality and the cup product,
without any recourse to the chain level.

Goresky and MacPherson returned to the chain-level intersection pairing as a
tool for constructing an intersection pairing in intersection homology
\cite[Section 2]{GM}.  They introduced the PL chain complex $C_*M$ (as the 
direct limit of simplicial chains under subdivision; see \fullref{sec2} 
) and defined the intersection pairing (which they denoted by $\cap$) on 
a certain subset of $C_*M\times C_*M$ by means of an elegant construction in 
which the procedure of the previous paragraph is reversed: the chain-level 
intersection pairing is derived from the relative versions of Poincar\'e 
duality and the cup product.  Their version of the chain-level intersection 
pairing is probably equivalent to Lefschetz's second definition. 

In order to prove (or even state) \fullref{1.1} it is necessary to extend
the domain of the intersection pairing from a subset of $C_*M\times C_*M$ to a
subset of $C_*M\otimes C_*M$.
The obvious way to do this would be to consider elements 
\[
\sum C_i\otimes D_i
\]
in which every pair $(C_i,D_i)$ is in the domain of the Goresky--MacPherson 
intersection pairing $\cap$ and to define the intersection pairing on such an
element to be 
\[
\sum C_i\cap D_i.
\]
But it is not at all clear that this is well defined, and it also is not clear
how to determine when an element of $C_*M\otimes C_*M$ has the required form
(which would make it difficult to show that the domain of the operation is a
full subcomplex).

The definition to be given in \fullref{sec7} resolves both of these issues
by defining the intersection pairing (up to a dimension shift) as 
the composite of the exterior product 
\[
\varepsilon\co C_*M\otimes C_*M \to C_*(M\times M)
\]
(see \fullref{sec6}) and the chain-level umkehr map 
\[
\Delta_!\co C^\Delta_*(M\times M)\to C_* M
\]
induced by the diagonal (see \fullref{sec5}); here 
$C^\Delta_*(M\times M)$ denotes the set of chains $E$ in $C_*(M\times M)$ for 
which both $E$ and $\partial E$ are in general position with respect to 
the diagonal.  With this definition, the domain of the intersection pairing 
(up to a dimension shift) is 
\[
\varepsilon^{-1}(C^\Delta_*(M\times M)).
\]
The analog of equation \eqref{1a.1} is immediate from the fact that
$\varepsilon$ and $\Delta_!$ are chain maps.

\section{PL chains}

\label{sec2}

We begin by reviewing some basic definitions.

A {\it simplicial complex\/} $K$ is a set of simplices in $\R^n$ (for some 
$n$) with two properties: every face of a simplex in $K$ is in $K$ and 
the intersection of two simplices in $K$ is a common face.  (A face of a 
simplex $\sigma$ is the simplex spanned by some subset of the vertices of 
$\sigma$.)

The {\it simplicial chain complex of $K$\/}, denoted $c_* K$, is defined by 
letting $c_p K$ be generated by pairs $(\sigma,o)$, where $\sigma$ is a 
$p$--simplex of $K$ and $o$ is an orientation of $\sigma$, subject to the 
relation $(\sigma,o)=-(\sigma,-o)$ where $-o$ denotes the opposite 
orientation.  We leave it as an exercise to formulate the definition of the 
boundary map $\partial$ (or see Spanier \cite[page 159]{Spanier}).  If we choose
orientations for the simplices of $K$ (with no requirement of consistency among
the orientations) then every nonzero element $c$ of $c_* K$ can be written 
uniquely in the form $\sum n_i\sigma_i$ with all $n_i\neq 0$.

The {\it realization\/} of $K$, denoted $|K|$, is the union of the simplices
of $K$.

A {\it subdivision\/} of $K$ is a simplicial complex $L$ with two properties:  
$|L|=|K|$ and every simplex of $L$ is contained in a simplex of $K$.  

The {\it subdivision category\/} of $K$ has an 
object for each subdivision $L$ of $K$ and a morphism $L\to L'$ whenever $L'$ 
is a subdivision of $L$. 

If $L'$ is a subdivision of $L$ there is an induced monomorphism 
$c_* L\to c_* L'$ which takes $(\sigma,o)$ to $\sum (\tau,o_\tau)$,
where the sum runs over all $\tau\in L'$ which are contained in $\sigma$ and
have the same dimension as $\sigma$, and $o_\tau$ is the orientation induced 
by $o$.  This makes $c_*$ a covariant functor on the subdivision category of
$K$.

A subspace $X$ of $\R^n$ will be called
a {\it PL space\/} if there is a simplicial complex $K$ with $X=|K|$.
$K$ will be called a {\it triangulation\/} of $X$;
note that $X$ determines $K$ up to subdivision by
Bryant \cite[page 222]{Bryant}.

The {\it PL chain complex\/} of a PL space $|K|$, denoted $C_*|K|$, is the 
direct limit
\[
\colim_L c_* L
\]
taken over the subdivision category of $K$. 

\begin{remark}
This definition is taken from Goresky and MacPherson \cite[Subsection 1.2]{GM}, which seems to be the
first place where the PL chain complex was defined.
\end{remark}

Note that the direct system defining $C_*|K|$ is a rather simple one: the 
subdivision category is a directed set (because any two subdivisions have a 
common refinement \cite[page 222]{Bryant}), and all of the maps $c_*L\to 
c_* L'$ are monomorphisms.  It follows that each of the maps $c_* L\to 
C_*|K|$ is a monomorphism.

\begin{remark}
The homology of $c_* L$ is canonically isomorphic to the singular homology of 
$|K|$ by \cite[Theorems 4.3.8 and 4.4.2]{Spanier}; since homology commutes with 
colimits over directed sets, the homology of $C_* |K|$ is also canonically 
isomorphic to the singular homology of $|K|$.
\end{remark}

Now let $C$ be a nonzero element of $C_* |K|$.  There is a subdivision $L$ 
of $K$ with $C$ in $c_* L$, so (after choosing orientations for the simplices 
in $L$) we can write $C=\sum n_i\sigma_i$ where the $\sigma_i$ are 
simplices in $L$ and the $n_i$ are nonzero.  We define the {\it support of 
$C$}, denoted $\supp(C)$, to be $\bigcup \,\sigma_i$; this is independent of 
the choice of $L$. The support of $0$ is defined to be the empty set.

\section{A useful lemma}

\label{sec3}

Let $K$ be a simplicial complex.  A {\it subcomplex\/} of $K$ is a subset 
$K'$ of $K$ with the property that every face of every simplex in $K'$ is also
in $K'$.

A {\it PL subspace\/} of $|K|$ is a space of the form $|L|$ where $L$ is a 
subcomplex of a subdivision of $K$.

The next lemma is taken from Section 1.2 of \cite{GM};  it gives a way of using
relative homology to make chain-level constructions.

\begin{lemma}
\label{3.1}
Let $K$ be a simplicial complex and let $A$ and $B$ be PL subspaces of $|K|$
such that $B\subset A$ and $\dim B\leq \dim A-1$.  Let $p=\dim A$. 

{\rm (a)}\qua There 
is a natural  isomorphism $\alpha_{A,B}$ from $H_p(A,B)$ to the abelian group
\[
\{\, C\in C_p(K) \,:\, 
\supp(C)\subset A \text{\rm\ and } \supp(\partial C)\subset B \,\}.
\]
{\rm (b)}\qua The following diagram commutes: 
\[
\xymatrix{
H_p(A,B)
\ar[d]_{\partial}
\ar[r]^-{\alpha_{A,B}}
&
\{\, C\in C_p(K) \,:\,
\supp(C)\subset A \text{\rm\ and } \supp(\partial C)\subset B \,\}
\ar[d]^{\partial} 
\\
H_{p-1} (B,\emptyset)
\ar[r]^-{\alpha_{B,\emptyset}}
&
\{\, D\in C_{p-1}(K) \,:\,
\supp(D)\subset B \text{\rm\ and } \partial D=0
\,\} 
}
\]
\end{lemma}

\begin{proof}
For part (a), note that
$H_p(A,B)$ is isomorphic to the $p$--th homology of the complex $C_*A/C_* B$,
and this in turn is isomorphic to the quotient of the relative cycles by the
relative boundaries.  The set specified in the lemma is the set of
relative cycles, while the set of relative boundaries is
$\partial(C_{p+1} A )+C_p B$, which is zero because of the hypotheses.
Part (b) is immediate from the definitions.
\end{proof}

\section{An umkehr map in relative homology} 

\label{sec4}

A {\it PL map\/} from $|K|$ to
$|K'|$ is a continuous function $f$ with the property that, for some
subdivision $L$ of $K$, the restriction of $f$ to each simplex of $L$ is
an affine map with image in a simplex of $K'$.  

A {\it PL homeomorphism\/} is a PL map which is a homeomorphism.

An {\it $m$--dimensional PL manifold\/} is a PL space $M$ with the property that 
each point of $M$ is contained in the interior of a PL subspace which is PL 
homeomorphic to the $m$ simplex. 

Let $M$ be a compact oriented $m$--dimensional PL manifold and let $A$ and 
$B$ be PL subspaces of $M$ with $B\subset A$.  Let $N$ be a compact oriented 
PL manifold of dimension $n$ and let $f\co N\to M$ be a PL map.  Let
$A'=f^{-1}(A)$ and $B'=f^{-1}(B)$.

We want to construct a map
\begin{equation}
\label{e4.1}
f_!\co H_*(A,B)\to H_{*+n-m}(A',B')
\end{equation}
(one should think of this as taking a homology class to its inverse image 
with respect to $f$).

%

Let $(U',V')$ be an open pair in $N$ with $A'\subset U'$ and 
$B'\subset V'$.  Choose an open pair $(U,V)$ in $M$ with 
$A\subset U$, $B\subset V$,
$f^{-1}(U)\subset U'$ and $f^{-1}(V)\subset V'$ (for example, we can let $U=M-f(N-U')$ and
$V=M-f(N-V')$).  Consider the composite
\begin{multline*}
H_*(A,B)\to
H_*(U,V)\cong
\Check{H}^{m-*} (M-V,M-U) 
\\
\stackrel{f^*}{\rightarrow}
\Check{H}^{m-*} (N-V',N-U') \cong
H_{*+n-m} (U',V'),
\end{multline*}
where the second and fourth maps are Poincar\'e--Lefschetz duality 
isomorphisms; 
see Dold \cite[Proposition VIII.7.2]{Dold}.
By the naturality of the cap product \cite[VIII.7.6]{Dold} this  
composite is independent of the choice of $(U,V)$ and is natural with respect 
to $(U',V')$.  We therefore get a map
\[
H_*(A,B)\to \lim H_{*+n-m} (U',V')
\]
where the inverse limit is taken over all open pairs 
$(U',V')\supset (A',B')$.  This inverse limit is 
isomorphic to $H_{*+n-m}(A',B')$ by 
\cite[Exercise 4 at the end of Section VIII.13]{Dold}; here we use the fact
that the realization of a simplicial complex is an ENR (see for example
\cite[Proposition IV.8.12]{Dold}).  This completes the construction of the map
\eqref{e4.1}.


For use in the next section, we need: 

\begin{lemma}
\label{4.1}
The following diagram commutes:
\[
\xymatrix{
H_*(A,B)
\ar[d]_{\partial}
\ar[r]^-{f_!}
&
H_{*+n-m} (A',B')
\ar[d]^{\partial}
\\
H_{*-1} (B)
\ar[r]^-{f_!}
&
H_{*+n-m-1} (B')
}
\]
\end{lemma}

\begin{proof}
This follows easily from \cite[VII.12.22]{Dold}. 
\end{proof}

\section{An umkehr map at the chain level}

\label{sec5}

Let $M$, $N$ and $f\co N\to M$ be as in the previous section.

We say that a PL subspace $A$ of $M$ is in 
{\it general position\/} with respect
to $f$ if
\[
\dim(f^{-1}(A))\leq \dim A +n -m.
\]
(The dimension of the empty set is defined to be $-\infty$, so if $f^{-1}(A)$
is empty then $A$ is in general position.)

\begin{remark}
\label{5.1}
For later use we make two observations.

(a)\qua Suppose that $f$ is a composite $gh$, that $A$ is in general position with
respect to $g$, and that $g^{-1}(A)$ is in general position with respect to
$h$.  Then $A$ is in general position with respect to $f$.

(b)\qua Suppose that $N$ is a Cartesian product $M\times M_1$ and $f\co N\to M$ 
is the 
projection.  Then every $A$ is in general position with respect to $f$.
\end{remark}

A $p$--chain $C$ in $C_* M$ is said to be in general position with
respect to $f$ if
\[
\dim(f^{-1}(\supp(C)))\leq p+n-m.
\]
Let $C^f_* M$ be the set of all chains $C\in C_* M$ for which both $C$ and
$\partial C$ are in general position with respect to $f$.  Note that $C^f_* 
M$ is a subcomplex of $C_* M$.

We want to construct a chain map
\[
f_!\co  C^f_* M\to C_{*+n-m} N.
\]
So let $C\in C^f_q M$.  Let $[C]$ be the homology class of $C$ in 
$H_q(\supp(C),\supp(\partial C))$.  
Let $T$ be the abelian group
\[
\{\, D\in C_{q+n-m} N \,|\, \supp(D)\subset f^{-1}(\supp(C)) \text{\rm\ and }
\supp(\partial D)\subset f^{-1}(\supp(\partial C)) \,\}.
\]
We define $f_!(C)$ to be the image of $[C]$ under the
following composite:
\begin{multline*}
H_q(\supp(C),\supp(\partial C)) \stackrel{f_!}{\rightarrow}
H_{q+n-m}(f^{-1}(\supp(C)),f^{-1}(\supp(\partial C)))
\\
\cong
T
\hookrightarrow
C_{q+n-m}N
\end{multline*}
Here the first map was constructed in \fullref{sec4} and the isomorphism 
is from \fullref{3.1} (which applies because of the hypothesis that both 
$C$ and $\partial C$ are 
in general position with respect to $f$).  $f_!$ is a chain map by 
\fullref{3.1}(b) and \fullref{4.1}.

\begin{remark}
\label{5.2}
Note that, by definition of $T$, we have 
$
\supp(f_!(C))\subset
f^{-1}(\supp(C))
$.
\end{remark}

\section{The exterior product for PL chains}

\label{sec6}

Let $\sigma_1$ and $\sigma_2$ be simplices.
It is easy to see that $\sigma_1\times \sigma_2$ is a PL space;
that is, there is a simplicial complex $J$ with $|J|=\sigma_1\times 
\sigma_2$.  Note that there is no canonical way to choose $J$, but that any 
two choices of $J$ have a common subdivision.

It follows that the product of any two PL spaces is a PL space.

Let $|K_1|$ and $|K_2|$ be PL spaces.
We want to construct a map
\begin{equation}
\label{e6.1}
\varepsilon\co C_* |K_1| \otimes C_* |K_2| \to C_*(|K_1|\times |K_2|),
\end{equation}
called the {\it exterior product\/}.

As a first step, let $L_1$ and $L_2$ be subdivisions of $K_1$ and $K_2$
respectively.  We will define a map
\begin{equation}
\label{e6.2}
\varepsilon'\co c_*L_1\otimes c_*L_2 \to C_*(|K_1|\times |K_2|)
\end{equation}
(see \fullref{sec2} for the definition of $c_*$).

It suffices to define $\varepsilon'$ on generators, so for $i=1,2$ let
$\sigma_i$ be a simplex of $L_i$ with orientation $o_i$.
Let $J$ be a simplicial complex with 
\[
|J|=\sigma_1\times\sigma_2.
\]
Then $\varepsilon'((\sigma_1,o_1)\otimes(\sigma_2,o_2))$ is defined to be
\[
\sum (\tau,o_\tau)
\]
where $\tau$ runs through the simplices of $J$ with dimension
$\dim \sigma_1 +\dim \sigma_2 $, and $o_\tau$ is the orientation of $\tau$ 
induced by $o_1\times o_2$.

The maps $\varepsilon'$ are consistent as $L_1$ and 
$L_2$ vary; passage to colimits gives the map $\varepsilon$.

\begin{remark}
\label{6.2}
(a)\qua It is easy to check that $\varepsilon$ is a monomorphism.

(b)\qua 
The quasi-isomorphism relating $c_*$ to singular chains 
\cite[Theorems 4.3.8 and 4.4.2]{Spanier}
takes $\varepsilon$ to the Eilenberg--MacLane shuffle product 
\cite[VI.12.26.2]{Dold}.  Since the latter is a quasi-isomorphism, so is
$\varepsilon$.

(c)\qua
For singular chains, the shuffle product has an explicit natural homotopy
inverse, namely the Alexander--Whitney map \cite[VI.12.26.2]{Dold}.
Unfortunately the Alexander--Whitney map requires an ordering of the vertices 
of each simplex, so it seems to have no analog for PL chains.
\end{remark}

\section{The chain-level intersection pairing}

\label{sec7}

We now have the ingredients needed to define the chain-level intersection
pairing.

Let $M$ be a compact oriented PL manifold of dimension $m$ and let 
$\Delta$ 
be the diagonal map from $M$ to $M\times M$.  As in \fullref{sec5}, let $C^\Delta_* 
(M\times M)$ be the subcomplex of $C_* (M\times M)$ consisting of chains $C$
for which both $C$ and $\partial C$
are in general position with respect to $\Delta$.

It is convenient to shift degrees so that the intersection pairing preserves
degree.  For a chain complex
$C_*$ and an integer $n$, we will write $\Sigma^n C_*$
for the $n$--fold suspension of $C_*$, that is, the chain complex with 
$C_{i}$ in degree $i+n$.
\[
G_2 \subset \Sigma^{-2m}(C_* M\otimes C_* M) \leqno{\hbox{Define}}
\]
to be $\Sigma^{-2m}(\varepsilon^{-1}(C^\Delta_* (M\times M)))$, where 
$\varepsilon$ is 
the exterior product (the $G$ stands for ``general position'' and the 
subscript 2 will be explained in \fullref{sec9}).

The chain-level intersection pairing $\mu$ is the composite
\[
G_2 \labarrow{\varepsilon}
\Sigma^{-2m}C^\Delta_*(M\times M) 
\labarrow{\Delta_!}
\Sigma^{-m}C_* M.
\]

\begin{remark}\label{rem8.1}
It is not difficult to check that, if $C$ and $D$ are chains for which the
Goresky--MacPherson intersection pairing $C\cap D$ is defined 
\cite[pages 141--142]{GM}, then (up to the dimension shifts in the
definitions of $G_2$ and $\mu$) $C\otimes D$ is in $G_2$ and 
$\mu(C\otimes D)=C\cap D$.
\end{remark}

\section{Leinster partial commutative DGAs}

\label{sec8}

Our main goal in the rest of the paper is to show that the 
chain-level intersection pairing and its iterates determine a partially 
defined commutative DGA structure on $\Sigma^{-m}C_*M$.

First we need a precise definition of ``partially defined commutative
DGA.'' We will use the definition given by Leinster in \cite[Section 
2.2]{Leinster} (but note that Leinster uses the term ``homotopy algebra'' 
instead of ``partial algebra'').  

\begin{notation}
\label{8.05}

(i)\qua For $k\geq 1$ let $\bar{k}$ denote the set $\{1,\ldots,k\}$.  Let
$\bar{0}$ be the empty set.  

(ii)\qua Let $\Phi$ be the full subcategory of $\mathrm{Set}$ with objects 
$\bar{k}$ for $k\geq 0$.  

(iii)\qua Given a functor $A$ defined on $\Phi$, write $A_k$ (instead of 
$A({\bar{k}})$) for the value of $A$ at ${\bar{k}}$.
\end{notation}

Disjoint union gives a functor $\coprod\co \Phi\times\Phi\to \Phi$.  
In particular, if $A$ is a functor defined on $\Phi$ then the functor 
$A\circ \coprod$ on $\Phi\times\Phi$ takes $(\bar{k},\bar{l})$ to 
$A_{k+l}$.

\begin{notation}
\label{8.06}

(i)\qua Let $\Ch$ denote the category of chain complexes.  

(ii)\qua Let $({\mathbb Z},0)$ denote the chain complex which has $\mathbb Z$ in 
dimension $0$ and $0$ in all other dimensions.
\end{notation}

\begin{definition}
\label{8.1}
A Leinster partial commutative DGA is a functor $A$ from $\Phi$ to $\Ch$ 
together with chain maps
\[
\xi_{k,l}\co A_{k+l}\to A_k\otimes A_l
\quad\text{for each $k,l$\qquad and}\qquad
\xi_0\co A_0\to ({\mathbb Z},0)
\]
such that the following conditions hold.

(i)\qua The collection $\xi_{k,l}$ is a natural transformation from 
$A\circ\coprod$ to $A\otimes A$ (considered as functors from $\Phi\times \Phi$ 
to $\Ch$).

(ii)\qua The following diagram commutes for all $k,l,n$:
\[
\xymatrix{
A_{k+l+n}
\ar[r]^-{\xi_{k+l,n}}
\ar[d]_{\xi_{k,l+n}}
&
A_{k+l}\otimes A_n
\ar[d]^{\xi_{k,l}\otimes 1}
\\
A_k\otimes A_{l+n}
\ar[r]^-{1\otimes \xi_{l,n}}
&
A_k\otimes A_l\otimes A_n
}
\]
(iii)\qua Let $\tau\co \overline{k+l}\to \overline{k+l}$ be the block permutation 
that transposes $\{1,\ldots,k\}$ and $\{k+1,\ldots,k+l\}$.  The following diagram
commutes for all $k,l$:
\[
\xymatrix{
A_{k+l}
\ar[r]^-{\xi_{k,l}}
\ar[d]_{\tau_*}
&
A_k\otimes A_l
\ar[d]^{\cong}
\\
A_{k+l}
\ar[r]^-{\xi_{l,k}}
&
A_l\otimes A_k
}
\]
(iv)\qua The following diagram commutes for all $k$:
\[
\xymatrix{
A_k
\ar[r]^-{\xi_{0,k}}
\ar[dr]_{\cong}
&
A_{{0}}\otimes A_k
\ar[d]^{\xi_0\otimes 1}
\\
&
{\mathbb Z}\otimes A_k
}
\]
(v)\qua $\xi_0$ and each $\xi_{k,l}$ are quasi-isomorphisms.
\end{definition}

\begin{remark} 
\label{8.3}
(a)\qua
An ordinary commutative DGA $B$ determines a Leinster partial commutative DGA
with $A_k=B^{\otimes k}$.

(b)\qua Conversely, the proof of \cite[Theorem 1]{Wilson} can be modified to show 
that Leinster partial commutative DGAs can be functorially replaced by 
quasi-isomorphic $E_\infty$ DGAs.  
\end{remark}

\begin{remark}
\label{rev2}
(a)\qua
\fullref{8.1} is the precise analog, for the category $\Ch$, of
Segal's $\Gamma$--spaces \cite{Segal}.  This is not immediately obvious,
since a $\Gamma$--space is a functor on the category $\F$ of {\it based\/} finite
sets \cite[Remark 1.4]{MT}; the point is that the maps $\xi_{k,l}$ in 
\fullref{8.1} encode the same information as the projection maps in 
Segal's definition.

(b)\qua It will be shown in \cite{McClure} that partial commutative DGAs in 
the sense of Kriz and May \cite[Section II.2]{KM} are the same thing as
Leinster partial commutative DGAs in which all of the maps $\xi_{k,l}$ are 
monomorphisms.  
\end{remark}

\section[The functor G]{The functor $G$}

\label{sec9}

As a first step in showing that the intersection pairing on 
$\Sigma^{-m}C_*M$ 
extends to a Leinster partial commutative DGA
structure, we define a functor $G$ from $\Phi$ to $\Ch$ with
$G_1=\Sigma^{-m}C_*M$.  The $G$ stands for ``general
position.''

$G_2$ has already been defined in \fullref{sec7}.
To define $G_k$ for $k\geq 3$ we need a definition.

Let $R\co \bar{k}\to\bar{k'}$ be any
map. Define
\[
R^*\co  M^{k'}\to M^{k}
\]
to be the composite
\[
M^{k'}=\Map(\bar{k'},M)\to
\Map(\bar{k},M)
=M^k
\]
where the second arrow is induced by $R$.  Thus
the projection of $R^*(x_1,\ldots,x_{k'})$ on the $i$--th factor
is $x_{R(i)}$. 

If $R\co \bar{k}\twoheadrightarrow\bar{k'}$ is a 
surjection then we think of $R^*$ as a generalized diagonal map.  For 
example, if $k'$ is $1$ and $R$ is the constant map then $R^*\co M\to M^{k}$ is 
the usual diagonal map.  

Let $\varepsilon_k$ denote the $k$--fold exterior product
\[
(C_*M)^{\otimes k}\hookrightarrow 
C_*(M^k).
\]

\begin{definition}\label{def10.1}
Define $G_0$ to be ${\mathbb Z}$ and $G_1$ to be $\Sigma^{-m}C_*M$.  For 
$k\geq 2$ 
define $G_k$ to be the subcomplex of $\Sigma^{-mk}((C_*M)^{\otimes k})$ 
consisting of the elements $\Sigma^{-mk}C$ for which both $\varepsilon_k(C)$ 
and $\varepsilon_k(\partial C)$ are in
general position with respect to all generalized diagonal maps, that is, 
\[
G_k=\bigcap_{k'<k}
\bigcap_{R\co \bar{k}\twoheadrightarrow\bar{k'}  }
\Sigma^{-mk}(\varepsilon_k^{-1}(C_*^{R^*}M^{k})).
\]
\end{definition}

\begin{remark}
One might at first expect a simpler definition of $G_k$, in which general 
position is required only with respect to the ordinary diagonal $M\to M^k$.
The more complicated definition given here is needed for \fullref{10.05}.
\end{remark}

\begin{lemma}
\label{9.05}
If $\Sigma^{-mk}C$ is in $G_k$ then both $\varepsilon_k(C)$ and 
$\varepsilon_k(\partial C)$ are in general position with respect 
to all maps $R^*$.
\end{lemma}

\begin{proof}
Any $R$ factors as $R_1R_2$, where $R_1$ is an inclusion and $R_2$ is a
surjection.  Then $R^*=R_2^*R_1^*$, and $R_1^*$ is a composite of
projection maps.  The lemma now follows from \fullref{5.1}(a) and 
\fullref{5.1}(b).
\end{proof}

It remains to define the effect of $G$ on morphisms in $\Phi$.  First 
we need three lemmas.

\begin{lemma}
\label{9.1}
Let 
$
N_1\labarrow{g} N_2\labarrow{f} N_3
$
be a diagram of compact oriented PL manifolds and PL maps.  Let 
$C\in C_*^f N_3\cap C_*^{fg} N_3$.
Then 
$${\rm (a)}\qua f_! C\in C_*^g N_2 \qquad\text{and} \qquad {\rm (b)}\qua  g_!f_!C=(fg)_! C.$$
\end{lemma}

\begin{proof}
Part (a) is immediate from \fullref{5.2} and (b) 
follows from the definitions. 
\end{proof}

\begin{lemma}
\label{9.2}
Let $f\co N_1\to M_1$ and $g\co N_2\to M_2$ be PL maps between compact 
oriented PL manifolds.  Then

{\rm (a)}\qua the exterior product
\[
\varepsilon\co C_*M_1\otimes C_*M_2 \to C_*(M_1\times M_2)
\]
takes $C^f_*M_1\otimes C^g_*M_2$ to $C^{f\times g}_*(M_1\times M_2)$, and

{\rm (b)}\qua the diagram
\[
\xymatrix{
C^f_p M_1\otimes C^g_q M_2
\ar[d]_{f_!\otimes g_!}
\ar[r]^-{\varepsilon}
&
C^{f\times g}_{p+q}(M_1\times M_2)
\ar[d]^{(f\times g)_!}
\\
C_{p+n_1-m_1}N_1\otimes C_{q+n_2-m_2} N_2
\ar[r]^-{\varepsilon}
&
C_{p+q+n_1+n_2-m_1-m_2}(N_1\times N_2)
}
\]
commutes for all $p$ and $q$, where $m_i$ (resp.\ $n_i$) is the dimension of 
$M_i$ (resp.\ $N_i$).
\end{lemma}

\begin{proof}
Part (a) is obvious from the definitions. Part (b) follows from 
\cite[VII.12.17]{Dold}.
\end{proof}

\begin{lemma}
Let $R\co \bar{k}\to\bar{k'}$ be any map.  Then 
$(R^*)_!\circ\varepsilon_k$ takes $G_k$ to $\varepsilon_{k'}(G_{k'})$.
\end{lemma}

\begin{proof}
We need to prove two things: that 
\begin{align*}
(R^*)_!(\varepsilon_k(G_k))&\subset
\Sigma^{-mk'}\varepsilon_{k'}(C_*(M)^{\otimes k'})\\
(R^*)_!(\varepsilon_k(G_k)) &\subset \Sigma^{-mk'}C^{S^*}_*(M^{k'}) \tag*{\hbox{and that}}
\end{align*} 
for all surjections $S\co \bar{k'}\to\bar{k''}$.  The first 
follows from \fullref{9.2} and the second from \fullref{9.05} and 
\fullref{9.1}(a).
\end{proof}

We can now define the effect of $G$ on morphisms by letting
\[
G_R\co G_k\to G_{k'}
\]
be $\varepsilon_{k'}^{-1}\circ(R^*)_!\circ\varepsilon_k$ (recall that
$\varepsilon_{k'}$ is a monomorphism).
\fullref{9.1}(b) implies that $G_{R\circ S}=G_R\circ G_S$ for all $R$ and
$S$.

\section[The maps xi_{k,l}]{The maps $\xi_{k,l}$}

\label{sec10}

In order to construct the maps
\[
\xi_{k,l}\co G_{k+l}\to G_k\otimes G_l
\]
it suffices to prove the following lemma.

\begin{lemma}
\label{10.05}
The inclusion
\[
G_{k+l}\hookrightarrow 
\Sigma^{-m(k+l)}(C_*M)^{\otimes (k+l)} 
\cong
\Sigma^{-mk}(C_*M)^{\otimes k}
\otimes
\Sigma^{-ml}(C_*M)^{\otimes l}
\]
has its image in $G_k\otimes G_l$.  
\end{lemma}

We can then define $\xi_{k,l}$ to be the inclusion
$G_{k+l}\hookrightarrow G_k\otimes G_l$.

In order to prove \fullref{10.05} we need a criterion for deciding when an element 
of $\Sigma^{-mk}(C_*M)^{\otimes k} \otimes \Sigma^{-ml}(C_*M)^{\otimes l}$
is in $G_k\otimes G_l$; we will build up to this in stages, culminating in
\fullref{10.1}.

Let $K$ be a triangulation of $M^k$. Choose
orientations for the simplices of $K$ (with no consistency required among the
choices).  Recall that (since orientations have been chosen) $c_p K$ is the 
free abelian group generated by the $p$--simplices of $K$.  Let 
$R\co \bar{k}\to \bar{k'}$ be any map and define
$c_p(K,R)$
to be the free abelian group generated by the $p$--simplices of $K$ that are
{\it not\/} in general position with respect to $R^*$.
Let 
\[
\Upsilon_p^{K,R}\co  c_p K \to c_p(K,R)
\]
be the homomorphism which is the identity on $c_p(K,R)$ and which
takes the $p$--simplices that are in general position with respect to $R^*$ to
0.  
Let
\[
\Psi_p^{K,R}\co c_p K \to c_p(K,R)\oplus c_{p-1}(K,R)
\]
be $\Upsilon_p^{K,R}+\Upsilon_{p-1}^{K,R}\circ \partial$.  

As an immediate consequence of the definitions, we have:

\begin{lemma}
{\rm (a)}\qua An element of $c_p K$ is in general position with respect to 
$R^*$ if and only if it is in the kernel of $\Upsilon_p^{K,R}$.

{\rm (b)}\qua An element of $c_p K$ is in $C_*^{R^*}(M^k)$ (that is, the element 
and its boundary are both in general position with respect to $R^*$)
if and only if it is in the kernel of $\Psi_p^{K,R}$.

{\rm (c)}\qua An element of $c_p K$ is in 
\[
\bigcap_{k'<k} \bigcap_{R\co \bar{k}\twoheadrightarrow\bar{k'}  }
C_*^{R^*}(M^{k})
\]
if and only if it is in the kernel of 
\[
\sum_R \Psi_p^{K,R}\co c_p K\to \bigoplus_R \,(c_p(K,R)\oplus c_{p-1}(K,R)).
\]
\end{lemma}

Next let $L$ be a triangulation of $M^l$.  We would like to
characterize the elements of $\Sigma^{-mk}c_p K\otimes \Sigma^{-ml}c_q L$ 
that are in $G_k\otimes G_l$. First we need some algebra.

\begin{lemma}
For exact sequences of abelian groups
\[
0\to A\to B\labarrow{g} C
\qquad
\text{and} 
\qquad
0\to D\to E\labarrow{h} F 
\]
with $C$ and $F$ torsion free, 
$A\otimes D$ is the kernel of 
\[
g\otimes 1 + 1\otimes h\co  B\otimes E \to (C\otimes E)\oplus(B\otimes F).
\]
\end{lemma}

\begin{proof}
We may assume without loss of generality that
$g$ and $h$ are surjective.  
Then the diagram
\[
\xymatrix{
&
0
\ar[d]
&
0
\ar[d]
&
0
\ar[d]
&
\\
0
\ar[r]
&
A\otimes D
\ar[d]
\ar[r]
&
B\otimes D
\ar[d]
\ar[r]
&
C\otimes D
\ar[d]
\ar[r]
&
0
\\
0
\ar[r]
&
A\otimes E
\ar[d]
\ar[r]
&
B\otimes E
\ar[d]
\ar[r]
&
C\otimes E
\ar[d]
\ar[r]
&
0
\\
0
\ar[r]
&
A\otimes F
\ar[d]
\ar[r]
&
B\otimes F
\ar[d]
\ar[r]
&
C\otimes F
\ar[d]
\ar[r]
&
0
\\
& 
0
&
0
&
0
&
}
\]
has exact rows and columns.  The lemma follows by an easy diagram chase.
\end{proof}

\begin{corollary}
\label{10.1}
Let $C\in (C_*M)^{\otimes k}\otimes (C_*M)^{\otimes l}$ be such that
$(\varepsilon_k\otimes\varepsilon_l)(C)$ is in $c_*K\otimes c_*L$.
Then $\Sigma^{-m(k+l)}C$ is in $G_k\otimes G_l$ if and only if 
$(\varepsilon_k\otimes\varepsilon_l)(C)$
is in the kernel of 
\[
\sum_{R\co \bar{k}\to \bar{k'}}
\Psi_p^{K,R}\otimes 1
+
\sum_{S\co \bar{l}\to \bar{l'}}
1\otimes \Psi_q^{L,S}
\]
\end{corollary}

\begin{proof}[Proof of \fullref{10.05}]
Let $\Sigma^{-m(k+l)}C\in G_{k+l}$. Then there are triangulations $K$ of 
$M^k$ and $L$ of
$M^l$ such that $(\varepsilon_k\otimes\varepsilon_l)(C)\in c_*K\otimes c_*L$.
Let $R\co \bar{k}\to \bar{k'}$.  Then both $\varepsilon_{k+l}(C)$ 
and $\varepsilon_{k+l}(\partial C)$ are
in general position with respect to $(R\times 1)^*$ (by definition of
$G_{k+l}$), and it is easy to see that this implies 
\begin{align*}
(\Psi_p^{K,R}\otimes 1)(\varepsilon_k\otimes\varepsilon_l)(C)&=0.\\
(1\otimes \Psi_q^{L,S})(\varepsilon_k\otimes\varepsilon_l)(C)&=0
\tag*{\hbox{Similarly}}
\end{align*}
for all $S$.  Thus $\Sigma^{-m(k+l)}C$ is in 
$G_k\otimes G_l$ by \fullref{10.1}.
\end{proof}

\section{The main theorem}

\label{sec11}

\begin{theorem}
\label{11.1}
The functor $G$ defined in \fullref{sec9} with the maps $\xi_{k,l}$ 
defined in \fullref{sec10} is a Leinster partial commutative DGA.
\end{theorem}

\begin{remark}
\label{rev18}
Since the maps $\xi_{k,l}$ are monomorphisms, $G$ is also a Kriz--May partial
commutative DGA (see \fullref{rev2}(b)).
\end{remark}

To prove \fullref{11.1}, we need to verify the five parts of \fullref{8.1}.  Part (i) follows easily from the definitions and \fullref{9.2}.
Parts (ii)--(iv) are immediate from the definition of $\xi_{k,l}$. 
Part (v) is an easy consequence of the following result,
which will be proved in the next section.

\begin{proposition}
\label{11.b.3}
The inclusion
\[
G_k\hookrightarrow \Sigma^{-mk} (C_*M)^{\otimes k}
\]
is a quasi-isomorphism for all $k$.  
\end{proposition}

\begin{remark}
\label{11.b.4}
When $k=2$ this is \fullref{1.1} of the introduction, up to the dimension
shift introduced in \fullref{sec7}.
\end{remark}

\section[Proof of Proposition 12.3]{Proof of \fullref{11.b.3}}

\label{sec12}

Throughout this section and the next we fix an integer $k\geq 2$.

A {\it PL homotopy\/} is a
PL map $h\co X\times I\to Y$, where $X$ and $Y$ are PL spaces and $I$ is 
the interval $[0,1]$ with its usual PL structure.

It will be convenient to have notation for the standard inclusion maps 
$X\to X\times I$.  We write $i_0$ (resp.\ $i_1$) for the map which takes 
$x$ to $(x,0)$ (resp.\ $(x,1)$).

We need a supply of PL homotopies that preserve the image of
\[
\varepsilon_k\co 
(C_* M)^{\otimes k}
\hookrightarrow 
C_*(M^k).
\]

\begin{definition}
Suppose that we are given a number $l$ with
$1\leq l\leq k$ and a PL homotopy
\[
\phi\co  M\times I\to M.
\]
The {\it $l$--th factor PL homotopy\/} determined by this data
is the composite
\[
M^k \times I
\cong
M^{l-1}\times (M\times I)\times M^{k-l}
\labarrow{1\times \phi\times 1}
M^k.
\]
\end{definition}


Let $\iota$ be the canonical element of $C_1(I)$. 

\begin{lemma}
\label{12.05}
Let 
$
h\co M^k\times I\to M^k
$
be an $l$--th factor PL homotopy for some $l$ and let $C$ be in the image of
$
\varepsilon_k\co 
(C_*M)^{\otimes k}
\hookrightarrow
C_*(M^k)
$.
Then 

{\rm (a)}\qua  $(h\circ i_1)_* (C)$ is in the image of $\varepsilon_k$, and

{\rm (b)}\qua $h_*(\varepsilon(C\otimes \iota))$ is in the image of 
$\varepsilon_k$.
\end{lemma}

\begin{proof}
This is an easy consequence of the definitions.
\end{proof}

For the proof of \fullref{11.b.3} we will use a filtration of 
$\Sigma^{-mk}(C_*M)^{\otimes k}$.

\begin{definition}
\label{12.1}
(i)\qua For $0\leq j\leq k$ define $\Lambda_j$ to be the set of all surjections
$R\co \bar{k}\twoheadrightarrow\bar{k'}$ 
such that for each $i>j$ the set $R^{-1}(R(i))$ has only one element.  

(ii)\qua For $0\leq j\leq k$ define $G_k^j$ to be 
the subcomplex of $\Sigma^{-mk}(C_*M)^{\otimes k}$
consisting of the chains $C$ for which both $\varepsilon_k(C)$ and
$\varepsilon_k(\partial C)$ are in
general position with respect to $R^*$ for all $R\in\Lambda_j$.
\end{definition}

Thus we have a filtration
\[
G_k=G_k^k\subset G_k^{k-1}\subset\cdots\subset
G_k^0=\Sigma^{-mk}(C_*M)^{\otimes k}.
\]
\fullref{11.b.3} follows immediately from:

\begin{proposition}
\label{12.b.05}
For each $1\leq j\leq k$ the inclusion $G_k^j\subset G_k^{j-1}$ is a 
quasi-iso\-mor\-phism.
\end{proposition}

For this we need a lemma which will be proved in \fullref{sec12a}.

\begin{lemma}
\label{12.b.1}
Suppose that $D\in \Sigma^{mk}G_k^{j-1}$ and $\partial D\in\Sigma^{mk}G_k^j$.  
Then there is a $j$--th factor homotopy 
\[
h\co  M^k\times I\to M^k
\]
such that 

{\rm (a)}\qua
$h\circ i_0$ is the identity,

{\rm (b)}\qua
the chains
$(h\circ i_1)_*(\varepsilon_k D)$, 
$(h\circ i_1)_*(\varepsilon_k (\partial D))$ and
$h_*(\varepsilon(\partial D\otimes \iota))$
are in general position with
respect to $R^*$ for all $R\in \Lambda_j$, and

{\rm (c)}\qua
the chain $h_*(\varepsilon(D\otimes \iota))$ is in general position 
with respect to $R^*$ for all $R\in \Lambda_{j-1}$.
\end{lemma}

\begin{proof}[Proof of \fullref{12.b.05}]
We have to show two things:

(i)\qua For $D$ a cycle in $\Sigma^{mk}G_k^{j-1}$,
there is a cycle $C$ in $\Sigma^{mk}G_k^j$ homologous to $D$.

(ii)\qua If $C$ is a cycle in $\Sigma^{mk}G_k^j$ which is the boundary of an
element of $\Sigma^{mk}G_k^{j-1}$ then $C$ is the boundary of
an element of $\Sigma^{mk}G_k^j$.

To show (i), choose a homotopy $h$ as in \fullref{12.b.1}.

\fullref{12.05} implies that $(h\circ i_1)_* (\varepsilon_k D)$ is in the 
image of $\varepsilon_k$, so we may define 
\[
C=\varepsilon_k^{-1}((h\circ i_1)_* (\varepsilon_k D)).
\]
$C$ is a cycle,
and \fullref{12.b.1}(b) implies that $C$ is in $\Sigma^{mk}G_k^j$.

\fullref{12.05} also implies that $h_*(\varepsilon(D\otimes
\iota))$ is in the image of $\varepsilon_k$, so we may define
\[
E=\varepsilon_k^{-1}(h_*(\varepsilon(D\otimes
\iota))).
\]
\fullref{12.b.1}(c) implies that $E$ is in $\Sigma^{mk}G_k^{j-1}$.

Let $\kappa,\lambda\in C_0 I$ be the $0$--chains associated to 
$0,1\in I$;  then 
$
\partial\iota=\lambda-\kappa
$.
Now
\begin{align*}
\varepsilon_k (\partial E)
&=\partial(h_*(\varepsilon(D\otimes \iota)))
\\
&=
h_*(\varepsilon(\partial D\otimes \iota))
+(-1)^{|D|} 
h_*(\varepsilon(D\otimes \lambda))
-(-1)^{|D|} 
h_*(\varepsilon(D\otimes \kappa))
\\
&= 
0+(-1)^{|D|} (h\circ i_1)_* (\varepsilon_k D)
-(-1)^{|D|} (h\circ i_0)_* (\varepsilon_k D)
\\
&=
(-1)^{|D|} \varepsilon_k(C-D).
\end{align*}
Since $\varepsilon_k$ is a monomorphism, this implies that $C$ is 
homologous to $D$.

To show (ii), let $D\in \Sigma^{mk}G_k^{j-1}$ with
$\partial D=C$.  Choose a homotopy $h$ as in \fullref{12.b.1}.
Then $(h\circ i_1)_*(\varepsilon_k D)$ and 
$h_*(\varepsilon(\partial D\otimes \iota))$
are in the image of $\varepsilon_k$ by \fullref{12.05}, so we may define
\[
E_1=\varepsilon_k^{-1}((h\circ i_1)_*(\varepsilon_k D))
\qquad
\text{and} 
\qquad
E_2=\varepsilon_k^{-1}(h_*(\varepsilon(\partial D\otimes \iota))).
\]
\fullref{12.b.1}(b) implies that $E_1$ and $E_2$ are in $\Sigma^{mk}G_k^j$.
Now
\begin{align*}
\varepsilon_k(\partial E_2)
&=
(-1)^{|D|+1}
(h_*(\varepsilon(\partial D \otimes \lambda))
-h_*(\varepsilon(\partial D \otimes \kappa)))
\\
&=
(-1)^{|D|+1}
((h\circ i_1)_*(\varepsilon_k\partial D)
-(h\circ i_0)_*(\varepsilon_k\partial D))
\\
&=
(-1)^{|D|+1}
\varepsilon_k
(\partial E_1
-C).
\end{align*}
Since $\varepsilon_k$ is a monomorphism, this implies
$
\partial((-1)^{|D|}E_2+E_1)
=C
$.
\end{proof}

\section[Proof of Lemma 13.5]{Proof of \fullref{12.b.1}}

\label{sec12a}

We will assume that $j=k$, since the other cases are essentially the 
same and the notation is simpler in this case.  So suppose we are given a $D$ 
satisfying:

\begin{assumption}
\label{14.04}
(i)\qua $D$ is in $\Sigma^{mk}G_k^{k-1}$.

(ii)\qua $\partial D$ is in $\Sigma^{mk}G_k^k$.  
\end{assumption}

With the assumption that $j=k$, \fullref{12.b.1} specializes to:

\begin{lemma} 
\label{14.01}
There is a $k$--th factor homotopy $h\co M^k\times I\to M^k$ such that

{\rm (a)}\qua
$h\circ i_0$ is the identity,

{\rm (b)}\qua
the chains
$(h\circ i_1)_*(\varepsilon_k D)$,
$(h\circ i_1)_*(\varepsilon_k (\partial D))$ and
$h_*(\varepsilon(\partial D\otimes \iota))$
are in general position with
respect to $R^*$ for all 
$R\co \bar{k}\twoheadrightarrow\bar{k'}$, 
and

{\rm (c)}\qua
the chain $h_*(\varepsilon(D\otimes \iota))$ is in general position
with respect to $R^*$ for all $R\in \Lambda_{k-1}$.
\end{lemma} 

\begin{remark}
\label{rev11}
Since $R^*$ is 1-1, the definition of general position simplifies
somewhat: a chain $C$ is in general position with respect to
$R^*$ if and only
if
\[
\dim(\supp(C)\cap\im(R^*))
\leq
\dim C+(k'-k)m.
\]
\end{remark}

Choose a triangulation $K$ of $M$ such that $D\in (c_*K)^{\otimes k}$. 

\begin{notation}
Let $\tau_1,\ldots,\tau_r$ be the simplices of $K$.  
\end{notation}

We fix orientations for $\tau_1,\ldots,\tau_r$ (with no requirement of 
consistency among the choices); this allows us to think of the $\tau_j$ as
generators of $c_*K$.

Since $D$ is in $(c_*K)^{\otimes k}$ it can be written as a sum 
\begin{equation}
\label{14.1}
D=\sum_\a n_\a \, \tau_{a_1}\otimes\cdots\otimes\tau_{a_k},
\end{equation}
where 
$\a$ runs through multi-indices 
$
(a_1,\ldots,a_k)\in \{1,\ldots,r\}^k
$
and $n_\a\in{\mathbb Z}$.  Then
\[
\supp(\varepsilon_k D)=\bigcup_{n_\a\neq 0} \,
\tau_{a_1}\times\cdots\times\tau_{a_k}.
\]
Similarly, $\partial D$ can be written as 
\begin{equation}
\label{14.2}
\partial D=\sum_\a n'_\a \, \tau_{a_1}\otimes\cdots\otimes\tau_{a_k},
\end{equation}
and we have
\begin{equation}
\label{14.3}
\supp(\varepsilon_k (\partial D))=\bigcup_{n'_\a\neq 0} \,
\tau_{a_1}\times\cdots\times\tau_{a_k}. 
\end{equation}
With this notation, we can spell out the meaning of \fullref{14.04}:

\begin{lemma}
\label{rev10}

{\rm (a)}\qua  If $S$ is a subset of $\overline{k-1}$ and $\a$ is a multi-index
with $n_\a\neq 0$ then
\[
\dim
\Bigl(\bigcap_{i\in S} \tau_{a_i}
\Bigr)
\leq
\sum_{i\in S} \dim \tau_{a_i}  + (1-|S|)m,
\]
where $|S|$ is the cardinality of $S$.

{\rm (b)}\qua The same inequality holds if $S$ is a subset of $\bar{k}$ 
and $\a$ is a multi-index with $n'_\a\neq 0$. 
\end{lemma}

\begin{proof}
For part (a), let $k'=k-|S|+1$ and let 
$R\co\bar{k}\twoheadrightarrow\ov{k'}$ be any surjection which takes $S$ to
$1$ and is 1-1 on the rest of $\ov{k}$.  
Note that 
$\supp(\varepsilon_k D)\cap \im(R^*)$ is homeomorphic to the subspace
\[
\bigcup_{n_\a\neq 0} \, \Bigl(\bigcap_{i\in S} \tau_{a_i} 
\times \prod_{i\not\in S} \tau_{a_i}
\Bigr)
\]
of $M^{k'}$.
By \fullref{14.04}(i) and \fullref{rev11} we have
\begin{equation}
\label{rev12}
\dim\Bigl(\bigcap_{i\in S} \tau_{a_i}
\times \prod_{i\not\in S} \tau_{a_i}
\Bigr)
\leq
\dim D + (1-|S|)m
\end{equation}
for all $\a$ with $n_\a\neq 0$.
But
\begin{gather*}
\dim\Bigl(\bigcap_{i\in S} \tau_{a_i}
\times \prod_{i\not\in S} \tau_{a_i}
\Bigr)
=
\dim\bigl(\bigcap_{i\in S} \tau_{a_i}\bigr)
+\sum_{i\not\in S}\, \dim \tau_{a_i},
\\
\tag*{\hbox{and}} 
\dim D
=
\sum_{i=1}^k \dim \tau_{a_i}.
\end{gather*}
Combining these equations
with inequality \eqref{rev12} completes the proof of part (a).  The proof of part (b) is 
similar.
\end{proof}

Next we need a general-position result that will be proved in
\fullref{sec13}--\fullref{sec13a}.

\begin{proposition}
\label{rev5}
Let $M$ be a compact PL manifold of dimension $m$ and let $K$ be a 
triangulation of $M$.  Then there is a PL homotopy 
\[
\phi\co M\times I \to M 
\]
with the following properties.

{\rm (a)}\qua $\phi\circ i_0$ is the identity.

{\rm (b)}\qua If $\sigma$ and $\tau$ are simplices of $K$ then 
$(\phi\circ i_1)(\sigma)$ and $\tau$ are in general position, ie
\[
\dim((\phi\circ i_1)(\sigma)\cap \tau)
\leq
\dim((\phi\circ i_1)(\sigma))
+
\dim \tau
-m.
\]
{\rm (c)}\qua 
If $\sigma$ and $\tau$ are any simplices of $K$ then 
\[
\dim(\phi(\sigma\times I)\cap \tau)\leq
\max( \dim\sigma+1+\dim\tau-m , \dim(\sigma\cap \tau)).
\]
\end{proposition}

\[
h\co M^k\times I\to M^k
\leqno{\hbox{We choose}}
\] 
to be the $k$--th factor homotopy determined by the
PL homotopy $\phi$ supplied by the Proposition.

Now let $R\co \bar{k}\twoheadrightarrow\bar{k'}$ be a surjection.  
To prove \fullref{14.01}
we need to show that the chains $(h\circ i_1)_*(\varepsilon_k D)$, 
$(h\circ i_1)_*(\varepsilon_k (\partial D))$, 
$h_*(\varepsilon(\partial D\otimes \iota))$, 
and (if $R\in \Lambda_{k-1}$) 
$h_*(\varepsilon(D\otimes \iota))$ 
are in general position with respect to $R^*$.  We will give the 
proof for $h_*(\varepsilon(\partial D\otimes \iota))$; the other cases are 
similar and easier.

Denote the set $R^{-1}(R(k))$ by $Q$.

First observe that $\supp(h_*(\varepsilon(\partial D\otimes \iota)))$
is the union of the $\dim D$ dimensional simplices of $h(\supp(\varepsilon_k
D)\times I)$.  In particular, it is contained in
the union over
all $\a$ such that
$n'_\a\neq 0$ of
\[
\tau_{a_1}\times\cdots\times\tau_{a_{k-1}}\times \phi(\tau_{a_k}\times I).
\]
It follows that $\supp(h_*(\varepsilon(\partial D\otimes \iota)))\cap
\im(R^*) $
is homeomorphic to a PL subspace of the union over
all $\a$ such that
$n'_\a\neq 0$ of
\[
\Bigl(
\prod_{j\neq R(k)}
\bigcap_{i\in R^{-1}(j)}
\tau_{a_i}
\Bigr)
\times
\Bigl(
\phi(\tau_{a_k}\times I)
\cap \bigcap_{i\in Q-\{k\}} \tau_{a_i}
\Bigr),
\]
and thus we have
\begin{multline*}
\dim(\supp(h_*(\varepsilon(\partial D\otimes \iota)))\cap
\im(R^*))
\\
\leq
\max_{n'_\a\neq 0}\,
\Bigl(
\sum_{j\neq R(k)} \,
\dim
\Bigl(
\bigcap_{i\in R^{-1}(j)}
\tau_{a_i}
\Bigr)
+
\dim
\Bigl(
\phi(\tau_{a_k}\times I)
\cap \bigcap_{i\in Q-\{k\}} \tau_{a_i}
\Bigr)
\Bigr).
\end{multline*}
It therefore suffices by \fullref{rev11} to show that for all $\a$ with $n'_\a\neq 0$,
\begin{multline}
\label{rev15}
\sum_{j\neq R(k)} \,
\dim
\Bigl(
\bigcap_{i\in R^{-1}(j)}
\tau_{a_i}
\Bigr)
+
\dim
\Bigl(
\phi(\tau_{a_k}\times I)
\cap \bigcap_{i\in Q-\{k\}} \tau_{a_i}
\Bigr)
\\
\leq
\dim(h_*(\varepsilon(\partial D\otimes \iota)))+(k'-k)m.
\end{multline}
Fix a multi-index $\a$ with $n'_\a\neq 0$.
By \fullref{rev5}(c) we know that one of the two following inequalities holds:
\begin{align}
\label{rev16}
\dim
\Bigl(
\phi(\tau_{a_k}\times I)
\cap \bigcap_{i\in Q-\{k\}} \tau_{a_i}
\Bigr)
&\leq
\dim \tau_{a_k} +1 +
\dim\Bigl( \bigcap_{i\in Q-\{k\}} \tau_{a_i} \Bigr) 
-m
\\
\label{rev17}
\dim
\Bigl(
\phi(\tau_{a_k}\times I)
\cap \bigcap_{i\in Q-\{k\}} \tau_{a_i}
\Bigr)
&\leq
\dim
\Bigl(
\tau_{a_k}
\cap \bigcap_{i\in Q-\{k\}} \tau_{a_i}
\Bigr)
\end{align}
If \eqref{rev17} holds then
\begin{align*}
\sum_{j\neq R(k)} \,
\dim
\Bigl(
&\bigcap_{i\in R^{-1}(j)}
\tau_{a_i}
\Bigr)
+
\dim
\Bigl(
\phi(\tau_{a_k}\times I)
\cap \bigcap_{i\in Q-\{k\}} \tau_{a_i}
\Bigr)
\\
&\leq
\sum_{j\neq R(k)} \,
\dim
\Bigl(
\bigcap_{i\in R^{-1}(j)}
\tau_{a_i}
\Bigr)
+
\dim
\Bigl(
\tau_{a_k}
\cap \bigcap_{i\in Q-\{k\}} \tau_{a_i}
\Bigr)
\\
&=
\sum_{j=1}^{k'} \,
\dim
\Bigl(
\bigcap_{i\in R^{-1}(j)}
\tau_{a_i}
\Bigr)
\\
&\leq
\dim(\supp(\varepsilon_k(\partial D))\cap
\im(R^*)) \quad\text{by equation (8)}
\\
&
\leq
\dim(\varepsilon_k(\partial D)) +(k'-k)m
\quad\text{by Assumption 14.1(ii)}
\\
&<
\dim(h_*(\varepsilon(\partial D\otimes \iota)))+(k'-k)m
\end{align*}
so inequality \eqref{rev15} holds in this case.

If \eqref{rev16} holds then we have
\begin{align*}
\sum_{j\neq R(k)} \,
&\dim
\Bigl(
\bigcap_{i\in R^{-1}(j)}
\tau_{a_i}
\Bigr)
+
\dim
\Bigl(
\phi(\tau_{a_k}\times I)
\cap \bigcap_{i\in Q-\{k\}} \tau_{a_i}
\Bigr)
\\
&\leq
\sum_{j\neq R(k)} \,
\Bigl(
\sum_{i\in R^{-1}(j)} \dim \tau_{a_i}  + (1-|R^{-1}(j)|)m
\Bigr)
\,
+
\dim(\phi(\tau_{a_k}\times I))
\\
&\qquad
+
\dim\Bigl( \bigcap_{i\in Q-\{k\}} \tau_{a_i} \Bigr)
-m
\quad\text{by \fullref{rev10}(b) and inequality \eqref{rev16}}
\\
&\leq
\Bigl(
\sum_{i\not\in R^{-1}(k)} \dim \tau_{a_i}
\Bigr)
+
(k'-1-k+|Q|)m 
+
\dim \tau_{a_k}+1
\\
&\qquad 
+
\dim\Bigl( \bigcap_{i\in Q-\{k\}} \tau_{a_i} \Bigr)
-m
\\
&\leq
\Bigl(
\sum_{i\not\in R^{-1}(k)} \dim \tau_{a_i}
\Bigr)
+
(k'-1-k+|Q|)m
+
\dim \tau_{a_k}+1
\\
&\qquad
+
\Bigl(
\sum_{i\in Q-\{k\}}\, \tau_{a_i} 
\Bigr)
+(1-|Q|+1)m
-m \quad\text{by \fullref{rev10}(b)}
\\
&=
\Bigl(
\sum_{i=1}^k \dim \tau_{a_i}
\Bigr)
+1
+(k'-k)m
\\
&=
\dim(\partial D) +1 +(k'-k)m
\quad\text{by equation \eqref{14.2}}
\\
&=
\dim(h_*(\varepsilon(\partial D\otimes \iota)))+(k'-k)m
\end{align*}
which proves inequality \eqref{rev15} in this case.

Thus we have shown that $h_*(\varepsilon(\partial D\otimes \iota))$
is in general position with respect to $R^*$.

\section[Background for the proof of Proposition 14.6]{Background for the proof of \fullref{rev5}}

\label{sec13}

First we have two simple facts about affine geometry which are the heart of 
the proof.  Recall that the {\it affine span\/} of a subset of $\R^n$ is
the smallest affine subspace containing it.

\begin{lemma}
\label{13.1}
Let $\sigma$ and $\tau$ be simplices in $\R^n$ 
such that the affine span of
$\sigma\cup \tau$ is all of $\R^n$.  Then $\sigma$ and $\tau$ are in general
position.
\end{lemma}

\begin{proof}
Let $U$ (resp.\ $V$) be the affine span of $\sigma$ (resp.\ $\tau$).
If $U\cap V$ is empty the statement is obvious.  Otherwise we can choose 
a point in $U\cap V$ and move it to the origin by a translation; then $U$ and 
$V$ become ordinary subspaces which span $\R^n$ and we have
\[
\dim(U\cap V)=\dim U+\dim V-n,
\]
which proves the lemma.
\end{proof}

\begin{notation}
If $\sigma$ is a simplex in $\R^n$ and $u$ is an element of $\R^n$ which is not
in $\sigma$,  the convex hull of 
$\sigma$ and $u$ will be denoted by $\langle \sigma,u\rangle$. 
\end{notation}

\begin{lemma}
\label{13.2}
Let $\sigma$ and $\tau$ be simplices in $\R^n$. 
Let $u$ be a point which is not in 
the affine span of $\sigma\cup\tau$. 
Then 
\[
\langle \sigma,u\rangle\cap \tau=\sigma\cap \tau.
\]
\end{lemma}

\begin{proof}
Let $v\in \langle\sigma,u\rangle\cap \tau$.  Since 
$v\in\langle\sigma,u\rangle$, we can write  $v$ in the form 
$\alpha u+(1-\alpha)s$, with $s\in\sigma$.  If $\alpha$ were nonzero we would
have 
\[
u=\frac{1}{\alpha}v-\frac{1-\alpha}{\alpha}s.
\]
Since $v\in\tau$, 
this would imply that $u$ is in the affine span of $\sigma\cup\tau$.
Therefore $\alpha$ must be 0, so $v$ is in $\sigma$, and 
hence in $\sigma\cap\tau$, which proves the lemma.
\end{proof}

Next we recall a well-known way of triangulating $\sigma\times I$.
By an {\it ordered\/} simplex we will mean a simplex with a total ordering of its
vertices.

\begin{lemma}
\label{13.3}
Let $\sigma$ be an ordered simplex and let $v_0<\cdots<v_l$ be the ordering of
its vertices.
For $0\leq i\leq l$ 
let $\sigma[i]\subset\sigma\times I$ be the convex hull of
\[
\{(v_j,0)\,|\, j\leq i\} \cup \{(v_j,1)\,|\, j\geq i\}.
\]
Then each $\sigma[i]$ is an $(l+1)$--simplex.

Also, 
the set $L$ whose elements are the $\sigma[i]$ and their faces
is a triangulation of $\sigma\times I$.
\end{lemma}

\begin{remark}
\label{13.4}
With the notation of \fullref{13.3}, let $\tau$ be the simplex spanned by
$v_0,\ldots,v_{l-1}$.  Then
\[
\sigma[i]=\langle\tau[i],(v_l,1)\rangle
\quad
\text{for each $i<l$,\qquad and}
\qquad
\sigma[l]=\langle \sigma\times \{0\}, (v_l,1)\rangle.
\]
\end{remark}

Finally, we need a tool for extending PL maps and homotopies.

\begin{construction}
\label{13.5}
Let $\rho$ be a simplex in $\R^n$ and let $u$ be an element 
of $\R^n$ which is not in $\rho$. 

Let $\Omega$ be a PL space with a PL homeomorphism $\omega\co \Omega\to 
\Delta^m$. 

(i)\qua
Let $f\co \rho\to \Omega$ be a PL map and $w$ an element of $\Omega$.
We can extend $\!f\!$ to a PL map
\[
\bar{f}\co \langle\rho,u\rangle\to \Omega
\]
by the formula
\[
\bar{f}(\alpha x+(1-\alpha) u)
=
\omega^{-1}(\alpha\omega(f(x))+(1-\alpha)\omega(w)).
\]
(ii)\qua
Next suppose we are given an ordering of the vertices of $\rho$; extend this
to $\langle\rho,u\rangle$ by letting $u$ be the maximal element.
Let $\phi\co \rho\times I\to \Omega$ be a PL homotopy 
and let $z$ and $z'$ be elements of $\Omega$.  
We can extend $\phi$ to a PL homotopy
\[
\bar{\phi}\co \langle\rho,u\rangle\times I\to \Omega
\]
as follows.  
Let $l-1$ be the dimension of $\rho$.  For $i<l$ we have 
\[
\langle \rho,u\rangle[i]=\langle\rho[i],(u,1)\rangle
\]
by \fullref{13.4}.  $\phi$ is already defined on $\rho[i]$, 
and we can extend it to $\langle\rho[i],(u,1)\rangle$ by using the
construction in part (i). For $i=l$, \fullref{13.4} gives 
\begin{equation*}
\langle \rho,u\rangle[l]
=
\langle \langle\rho,u\rangle\times \{0\}, (u,1)\rangle
=
\langle \langle\rho\times \{0\},(u,0)\rangle, (u,1)\rangle.
\end{equation*}
$\phi$ is already defined on $\rho\times \{0\}$, and we 
can extend it to 
$\langle \langle\rho\times \{0\},(u,0)\rangle, (u,1)\rangle$
by applying the
construction in part (i) twice, taking
$(u,0)$ to $z$ and $(u,1)$ to $z'$.
\end{construction}

\section[Proof of Proposition 14.6]{Proof of \fullref{rev5}}

\label{sec13a}

By the definition of PL manifold, there is a collection of PL subspaces 
$\Omega_i\subset M$ such that each $\Omega_i$ is PL homeomorphic to $\Delta^m$
and the interiors $\int(\Omega_i)$ cover $M$.
Choose PL homeomorphisms 
\[
\omega_i\co \Omega_i\to \Delta^m.
\]
Recall that by the definition in \fullref{sec2}, a PL space is 
given as a subspace of some $\R^n$, and therefore inherits a metric.  In 
particular, this is true for the PL manifold $M$.  Let us denote the metric 
on $M$ by $d$ and the standard norm on $\R^m$ by $||\ ||$.

\begin{definition}
\label{14.4}
(i)\qua For each $\Omega_i$, choose numbers $\gamma_i$ and $\delta_i$ with 
\[
||\omega_i(x)-\omega_i(y)||
\leq \gamma_i \, d(x,y)
\qquad
\text{and}
\qquad
d(x,y)
\leq \delta_i \, ||\omega_i(x)-\omega_i(y)||
\]
for all $x,y\in\Omega_i$ (such numbers exist because $\omega_i$ and its 
inverse are PL maps).

(ii)\qua Let $\lambda$ be the greater of $\max_i \gamma_i\delta_i$ and 1.
\end{definition}

\begin{definition}
\label{14.05}
Let $\eta$ be the Lebesgue number of the covering $\{\Omega_i\}$ (with 
respect to the metric $d$).
\end{definition}

Next observe that if \fullref{rev5} holds for some subdivision of $K$
then it holds for $K$. Choose a subdivision $L$ of $K$ such that

\begin{enumerate}
\renewcommand{\labelenumi}{(\roman{enumi})}
\item each $\Omega_i$ is a union of simplices of $L$,

\item the restriction of each $\omega_i$ to each simplex of
$L$ in $\Omega_i$ is affine,

\item the diameter of each simplex of $L$ is less than
$\frac{\eta}{2}$.
\end{enumerate}

It suffices to prove \fullref{rev5} for the triangulation $L$.

Choose an ordering 
$
v_1,\ldots,v_s
$
for the vertices of $L$.

\begin{definition}
\label{14.35}
For $1\leq p\leq s$ let $A_p$ be the union of the simplices of $L$ whose 
vertices are in the set $\{v_1,\ldots,v_p\}$.  Let $A_0$ be the empty set.
\end{definition}

Note that $A_s$ is $M$.

We will construct, by induction over $p$ with $0\leq p\leq s$, a PL homotopy
\[
\phi_p\co A_p\times I\to M
\]
with the following properties: 

\begin{enumerate}
\item
The restriction of $\phi_p$ to $A_{p-1}\times I$ is $\phi_{p-1}$.

\item
$\phi_p\circ i_0$ is the inclusion map of $A_p$ into $M$.

\item
For each $x\in A_p, t\in I$ we have
\[
d(\phi_p(x,t),x)\leq \frac{\eta}{2\lambda^{s-p}}.
\]
\item
If $\sigma$ is a simplex of $L$ in $A_p$ and $\tau$ is any simplex of $L$ then
$\phi(\sigma\times \{1\})$ and $\tau$ are in general position, and 
\[
\dim(\phi_p(\sigma\times I)\cap \tau)\leq
\max( \dim\sigma+1+\dim\tau-m , \dim(\sigma\cap \tau)).
\]
\end{enumerate}

This will complete the proof of \fullref{rev5}, because the homotopy 
$\phi_s$ will have the required properties.

The first step of the induction (the case $p=0$) is trivial. Suppose that 
$\phi_{p-1}$ has been constructed.  

\begin{notation}
\label{14.45}
Denote the simplices of $L$ which are in $A_p$ but not $A_{p-1}$ by 
\[
\pi_1,\ldots,\pi_t.
\]
For each 
$\pi_j$, let $\rho_j$ be the face opposite $v_p$; thus
\[
\pi_j=\langle \rho_j,v_p\rangle.
\]
\end{notation}

Combining property (iii) of the triangulation $L$ with property (3) of
$\phi_{p-1}$ and the fact that $\lambda\geq 1$, we see that for each $j$ the 
diameter of the set
\[
\pi_j\cup\phi_{p-1}(\rho_j\times I)
\]
is less than $\eta$.
It follows that for each $j$ we can choose a number $i(j)$ with 
\begin{equation}
\label{14.46}
\pi_j\cup\phi_{p-1}(\rho_j\times I) \subset \int(\Omega_{i(j)}).
\end{equation}

\begin{notation}
Let $\Xi$ denote the intersection of the sets $\int(\Omega_{i(j)})$.
\end{notation}

Note that $\Xi$ is nonempty (for example, it contains $v_p$).

If $z'$ is any point in $\Xi$
we can apply \fullref{13.5}(ii) 
(with $\Omega=\Omega_{i(j)}$ and $z=v_p$) to extend $\phi_{p-1}$ over all
$\pi_j\times I$ simultaneously. The resulting homotopy $\phi_p$ will 
automatically satisfy properties (1) and (2) above. We next state two lemmas
which will show that there is a $z'$ for which properties (3) and (4) hold.

\begin{lemma}
\label{rev3}
$\phi_p$ satisfies property {\rm (3)} if $z'$ is in the open ball 
$B$ of radius $\frac{\eta}{2\lambda^{s-p+1}}$ around $v_p$.
\end{lemma}

In order to verify property (4) for all simplices $\sigma$ of $L$ which are in
$A_p$, it suffices to consider the simplices which are in $A_p$ but not in
$A_{p-1}$ (that is, the simplices $\pi_1,\ldots,\pi_t$)
since the inductive hypothesis ensures that property (4) holds for all
simplices of $A_{p-1}$.

\begin{lemma}
\label{rev4}
For each $\pi_j$ and for each simplex $\tau$ of $L$, there is an open set
$U_{j,\tau}$ which is dense in $\Xi$ such that if $z'$ is in $U_{j,\tau}$ 
then

{\rm (a)}\qua
$\phi_p(\pi_j\times \{1\})$ and $\tau$ are in general position, and

{\rm (b)}\qua
$ \dim(\phi_p(\pi_j\times I)\cap \tau)\leq \max( \dim\pi_j+1+\dim\tau-m , 
\dim(\pi_j\cap \tau))$.
\end{lemma}

Before proving \fullref{rev3} and \fullref{rev4} we observe that
the set 
\[
U=B\cap \bigcap_{j,\tau} \,U_{j,\tau}
\]
will be dense in $B\cap \Xi$ (and in particular nonempty), and
if $z'$ is in $U$ then $\phi_p$ will satisfy properties (1)--(4), which
completes the inductive step and thereby the proof of \fullref{rev5}.

\begin{proof}[Proof of \fullref{rev3}]
Let $x\in \pi_j$ and $t\in I$.  Let $l$ be the
dimension of $\pi_j$. With the
notation of \fullref{13.3}, we have $(x,t)\in \pi_j[e]$ for some $e$ with
$0\leq e\leq l$.  There are two cases to consider: $e<l$ and $e=l$.  

In the first case, \fullref{13.4} allows us to write $(x,t)$ as
\[
\alpha \,(y,t')+(1-\alpha)\,(v_p,1)
\]
with $0\leq\alpha\leq 1$, $y\in \rho_j$ and $t'\in I$. By
\fullref{13.5}(ii) we have
\begin{equation}
\label{14.5}
\omega_{i(j)} (\phi_p(x,t))
=
\alpha\,\omega_{i(j)}(\phi_{p-1}(y,t')) +(1-\alpha)\,\omega_{i(j)}(z'),
\end{equation}
and by property (ii) of the triangulation $L$ we have
\begin{equation}
\label{14.6}
\omega_{i(j)}(x)=\alpha\,\omega_{i(j)}(y)+(1-\alpha)\,\omega_{i(j)}(v_p).
\end{equation}
Now we have
\begin{align*}
d(\phi_p(x,t),x)
&\leq
\delta_{i(j)}\,||\omega_{i(j)} (\phi_p(x,t))-\omega_{i(j)} (x) ||
\quad\text{by \fullref{14.4}(i)}
\\
&\leq
\delta_{i(j)}(\alpha\,||\omega_{i(j)}(\phi_{p-1}(y,t'))-\omega_{i(j)}(y)||
\\
&\qquad
+
(1-\alpha)\,||\omega_{i(j)}(z')-\omega_{i(j)}(v_p)||)
\\
&
\qquad\text{by equations \eqref{14.5} and \eqref{14.6}}
\\
&\leq
\delta_{i(j)}\gamma_{i(j)}
(\alpha\,d(\phi_{p-1}(y,t'),y)
+
(1-\alpha)\,d(z',v_p))
\\
&
\qquad\text{by \fullref{14.4}(i)}
\\
&\leq
\lambda
\left(\alpha\frac{\eta}{2\lambda^{s-p+1}}
+
(1-\alpha)\,d(z',v_p)\right)
\\
&
\qquad\text{by \fullref{14.4}(ii) and property (3) of $\phi_{p-1}$}
\end{align*}
and this will be $\leq \frac{\eta}{2\lambda^{s-p}}$
if $d(z',v_p)\leq \frac{\eta}{2\lambda^{s-p+1}}$.

For the second case of property (3) we have $e=l$, and \fullref{13.4} gives
\[
(x,t)=\alpha\,(y,0)+(1-\alpha)\,(v_p,1)
\]
for some $y\in\pi_j$. Equation
\eqref{14.6} is still valid, and equation \eqref{14.5} is replaced by 
\begin{equation}
\label{14.7}
\omega_{i(j)} (\phi_p(x,t))
=
\alpha\,\omega_{i(j)}(y) +(1-\alpha)\,\omega_{i(j)}(z').
\end{equation}
Now we have
\begin{align*}
d(\phi_p(x,t),x)
&\leq
\delta_{i(j)}\,||\omega_{i(j)} (\phi_p(x,t))-\omega_{i(j)} (x) ||
\quad\text{by \fullref{14.4}(i)}
\\
&\leq
\delta_{i(j)} (1-\alpha)\,||\omega_{i(j)}(z')-\omega_{i(j)}(v_p)||
\quad\text{by equations \eqref{14.6} and \eqref{14.7}}
\\
&\leq
\delta_{i(j)}\gamma_{i(j)} d(z',v_p)
\quad\text{by \fullref{14.4}(i)}
\\
&\leq
\lambda d(z',v_p)
\quad\text{by \fullref{14.4}(ii)} 
\end{align*}
For this to be $\leq \frac{\eta}{2\lambda^{s-p}}$
it again suffices to have $d(z',v_p)\leq \frac{\eta}{2\lambda^{s-p+1}}$.
\end{proof}

\begin{proof}[Proof of \fullref{rev4}]
We begin by considering the condition in part (a).
For this we need a precise description of the subspace 
$\phi_p(\pi_j\times\{1\})$. Recall that $\rho_j$ denotes the face of 
$\pi_j$ opposite to the vertex $v_p$.  By statement \eqref{14.46}, the set 
$\phi_{p-1}(\rho_j\times \{1\})$ is contained in $\Omega_{i(j)}$. The image of 
$\phi_{p-1}(\rho_j\times \{1\})$  under the PL homeomorphism
\[
\omega_{i(j)}\co \Omega_{i(j)}\to \Delta^m
\]
is a union of simplices which we will denote by $\chi_1,\ldots,\chi_u$.  
By \fullref{13.5}(ii), the 
subspace $\phi_p(\pi_j\times\{1\})$ is $\omega_{i(j)}^{-1}$ of the union of 
the simplices
\[
\langle \chi_1,\omega_{i(j)}(z') \rangle,
\ldots,
\langle \chi_u,\omega_{i(j)}(z') \rangle.
\]
In order for $z'$ to satisfy the condition in \fullref{rev4}(a), 
each of the pairs
\[
(\langle \chi_q,\omega_{i(j)}(z') \rangle,\omega_{i(j)}(\tau))
\]
must be in general position (note that $\omega_{i(j)}(\tau)$ is a simplex by
property (ii) of the triangulation $L$). By \fullref{13.1}, this condition is
automatically satisfied (with no restriction on $z'$) by those pairs for
which the affine span of $\chi_q\cup\omega_{i(j)}(\tau)$ is all of
$\R^m$.  For the remaining pairs, it suffices by \fullref{13.2} that
$\smash{\omega_{i(j)}(z')}$ should not be in the affine span of
$\chi_q\cup\omega_{i(j)}(\tau)$ (note that $\chi_q$ and
$\omega_{i(j)}(\tau)$ are in general position because we have assumed that 
$\phi_{p-1}$ satisfies property (4)).
Since this affine span is nowhere dense, the set of allowable $z'$ for each 
such pair is an open set $V_q$ which is dense in $\Xi$.  The
intersection of the $V_q$ is an open set $V$ which is dense 
in $\Xi$ and if $z'\in V$ then the condition in \fullref{rev4}(a) is 
satisfied.  

For part (b), let $l$ denote the dimension of $\pi_j$.  With the notation of 
\fullref{13.3} and \fullref{13.4}, we have 
\[
\pi_j\times I
=
\bigl(
\bigcup_{0\leq e < l} \langle\rho_j[e],(v_p,1)\rangle
\bigr)
\cup
\langle\pi_j\times\{0\},(v_p,1)\rangle.
\]
For each $e<l$ the image of $\phi_{p-1}(\rho_j[e])$  under $\omega_{i(j)}$
is a union of simplices which we will denote by $\psi^e_1,\ldots$. 
By \fullref{13.5}(ii), 
$\phi_p(\langle\rho_j[e],(v_p,1)\rangle)$ is 
$\omega_{i(j)}^{-1}$ of the union of the simplices
\[
\langle \psi^e_1,\omega_{i(j)}(z') \rangle,
\ldots
\]
and $\phi_p(\langle\pi_j\times\{0\},(v_p,1)\rangle)$
is $\omega_{i(j)}^{-1}$ of the simplex 
\[
\langle \omega_{i(j)}(\pi_j),\omega_{i(j)}(z')\rangle.
\]
In order for $z'$ to satisfy the condition in \fullref{rev4}(b), we must have
\begin{equation}
\label{rev6}
\dim(\langle \psi^e_q,\omega_{i(j)}(z') \rangle\cap 
\omega_{i(j)}(\tau))
\leq
\max( \dim\pi_j+1+\dim\tau-m , \dim(\pi_j\cap \tau))
\end{equation}
for all $q$, and 
\begin{multline}
\label{rev7}
\dim(\langle \omega_{i(j)}(\pi_j),\omega_{i(j)}(z')\rangle \cap
\omega_{i(j)}(\tau))
\\
\leq
\max( \dim\pi_j+1+\dim\tau-m , \dim(\pi_j\cap \tau)).
\end{multline}
Since we have assumed that $\phi_{p-1}$ satisfies property (4), we have
that for all $q$,
\begin{equation}
\label{rev8}
\dim(\psi^e_q\cap
\omega_{i(j)}(\tau))\leq
\max( \dim\rho_j+1+\dim\tau-m , \dim(\rho_j\cap \tau)).
\end{equation}
To prove inequality \eqref{rev6} we must consider two cases: either the affine
span of the union $\psi^e_q\,\cup \,\omega_{i(j)}(\tau)$ is all of $\R^m$ or it is not.  
In the first case we have (with no restriction on $z'$)
\begin{align*}
\dim(\langle \psi^e_q,\omega_{i(j)}(z') \rangle\cap
\omega_{i(j)}(\tau))
&\leq
\dim(\langle \psi^e_q,\omega_{i(j)}(z') \rangle) +\dim(\omega_{i(j)}(\tau))
-m
\\
&
\qquad\text{by \fullref{13.1}}
\\
&\leq
\dim \psi^e_q +1+\dim\tau -m
\\
&\leq
\dim \rho_j +2+\dim\tau -m
\\
&=
\dim \pi_j +1+\dim\tau -m.
\end{align*}
In the second case we assume that $\omega_{i(j)}(z')$ is not in the affine 
span of $\psi^e_q\,\cup \,\omega_{i(j)}(\tau)$.  Then we have
\begin{align*}
\dim(\langle \psi^e_q,\omega_{i(j)}(z') \rangle\cap
\omega_{i(j)}(\tau))
&=
\dim(\psi^e_q\cap \omega_{i(j)}(\tau))
\quad\text{by \fullref{13.2}}
\\
&\leq
\max( \dim\rho_j+1+\dim\tau-m , \dim(\rho_j\cap \tau))
\\
& \qquad\text{by inequality \eqref{rev8}}
\\
&\leq
\max( \dim\pi_j+1+\dim\tau-m , \dim(\pi_j\cap \tau)).
\end{align*}
To prove inequality \eqref{rev7} we again consider two cases: either the affine
span of the union $\omega_{i(j)}(\pi_j)\,\cup \, \omega_{i(j)}(\tau)$ is all of $\R^m$
or it is not.  In the first case we have (with no restriction on $z'$)
\begin{align*}
\dim(\langle \omega_{i(j)}(\pi_j),
&\omega_{i(j)}(z')\rangle \cap
\omega_{i(j)}(\tau))
\\
&\leq
\dim(\langle \omega_{i(j)}(\pi_j),\omega_{i(j)}(z') \rangle) 
+\dim(\omega_{i(j)}(\tau)) -m
\\
&
\qquad\text{by \fullref{13.1}}
\\
&\leq
\dim \pi_j +1+\dim\tau -m.
\end{align*}
In the second case we assume that $\omega_{i(j)}(z')$ is not in the affine 
span of the union $\omega_{i(j)}(\pi_j)\,\cup \, \omega_{i(j)}(\tau)$.  Then we have
\begin{align*}
\dim(\langle \omega_{i(j)}(\pi_j),\omega_{i(j)}(z')\rangle \cap
\omega_{i(j)}(\tau))
&=
\dim(\omega_{i(j)}(\pi_j)\cap
\omega_{i(j)}(\tau))
\\
&\qquad\text{by \fullref{13.2}}
\\
&=
\dim (\pi_j \cap \tau).
\end{align*}
To sum up, there is a dense open subset $W$ of $\Xi$ such that if $z'\in W$
then the condition of \fullref{rev4}(b) is satisfied.

Finally, let $U_{j,\tau}$ be $V\cap W$.
\end{proof}

{\small\parskip 0pt\vskip11pt minus 5pt\relax
{\sl \def\\{\futurelet\next\nocommawithnl}\def\nocommawithnl
  {\ifx\next\newline\else\unskip,\space\ignorespaces\fi}
  \theaddress\par}
{\rightskip0pt plus .4\hsize
{\def\tempab{}\tt\def~{\lower3.5pt\hbox{\char'176}}\def\_{\char'137}%
\ifx\theemail\tempab\else
  \vskip5pt minus 3pt\theemail\par\fi
  \ifx\theurl\tempab\else
  \vskip5pt minus 3pt\theurl\par\fi}
  \vskip11pt minus 5pt
            Proposed:\qua\theproposer\hfill
            Received:\qua\receiveddate\break
            Seconded:\qua\theseconders\hfill
        \ifx\reviseddate\tempab
              Accepted:\qua\accepteddate\else Revised:\qua
  \reviseddate\fi\break}}

\newpage


\title[Erratum for: On the chain-level intersection pairing for PL
manifolds]{Erratum for the paper\\`On the chain-level intersection pairing for PL
manifolds'}

\volumenumber{13}
\issuenumber{3}
\publicationyear{2009}
\papernumber{42}
\startpage{1775}
\endpage{1777}

\lognumber{0508}
\count0=1775

\doi{}
\MR{}
\Zbl{}

\received{18 July 2007}
\revised{}
\accepted{}
\published{12 March 2009}
\publishedonline{12 March 2009}
\proposed{}
\seconded{}
\corresponding{}
\editor{CPR}
\version{}



\renewcommand{\O}{{\cal O}}
\newcommand{\barm}{{\bar{\mu}}}
\newcommand{\Sm}{{\Sigma^{-m}}}
\newcommand{\Sn}{{\Sigma^{-n}}}


\begin{abstract}
Greg Friedman has pointed out that there are sign errors in the main
paper (above), and in 
particular \fullref{9.2}(b) (which is a key step in the proof of the main theorem)
is not correct as stated.

The purpose of this erratum is to provide a correction.
\end{abstract}

\maketitle
\section*{}\vglue -50pt
\addcontentsline{toc}{section}{Erratum}
\hypertarget{Err}{In} 
\fullref{sec4}, the umkehr map 
\[
H_p(A,B)\to H_{p+n-m}(A',B')
\]
should have a sign $(-1)^{(m-p)(n-m)}$ as in \cite[pages 314--315]{D} (see
\cite{F} for an explanation of where this sign comes from).  Also, it's
convenient to let the symbol $f_!$ stand for the desuspension
\[
\Sm H_*(A,B)\to \Sn H_*(A',B').
\]
Note that this map preserves degrees.

With these changes, \fullref{4.1} says that the diagram
\[
\xymatrix{
\Sm H_*(A,B)
\ar[r]^{f_!}
\ar[d]_{\partial}
&
\Sn H_*(A',B')
\ar[d]^{\partial}
\\
\Sm H_*(B)
\ar[r]^{f_!}
&
\Sn H_*(B')
}
\]
commutes.

In \fullref{sec7}, observe that if $C_*$ and $D_*$ are chain complexes and 
$m,n\in \mathbb Z$ then 
$(\Sigma^{-m} C_*)\otimes (\Sn D_*)$ 
and 
$\Sigma^{-(m+n)} (C_*\otimes D_*)$ 
will be {\it different} chain complexes 
when $n$ is odd: they are the same as graded abelian groups but have 
different differentials.   However, there is an isomorphism 
\[
\Theta:
(\Sigma^{-m} C_*)\otimes (\Sn D_*)
\to
\Sigma^{-(m+n)} (C_*\otimes D_*)
\]
which takes $\Sigma^{-m} x\otimes \Sn y$ to
$(-1)^{n|x|}\Sigma^{-(m+n)}(x\otimes y)$, and similarly for any number of 
tensor factors.

Given manifolds $M_1$ and $M_2$ of dimensions $m_1,m_2$, define
\[
\bar\varepsilon:(\Sigma^{-m_1} C_*M_1)\otimes (\Sigma^{-m_2} C_*M_2)
\to
\Sigma^{-(m_1+m_2)}C_*(M_1\times M_2)
\]
to be 
\[
(-1)^{m_1m_2}(\Sigma^{-(m_1+m_2)}\varepsilon)\circ\Theta.
\]
The motivation for this is that $\bar\varepsilon$
is Poincar\'e dual to the exterior product in cohomology (see \cite{F} for
details).  Similarly, given $M_1,\ldots,M_k$ define
\[
\bar\varepsilon:
\bigotimes 
\Sigma^{-m_i} C_*M_i
\to
\Sigma^{-\sum m_i} C_*(\textstyle{\prod} M_i)
\]
to be 
\[
(-1)^{e_2(m_1,\ldots,m_k)}(\Sigma^{-\sum m_i}\varepsilon)\circ\Theta
\]
where $e_2$ is the second elementary symmetric function (so that, for example,\break
$e_2(m_1,m_2,m_3)=m_1m_2+m_1m_3+m_2m_3$).

Now define 
\[
G_2\subset (\Sigma^{-m} C_*M)\otimes (\Sigma^{-m} C_*M)
\]
to be $\bar{\varepsilon}^{-1}(\Sigma^{-2m}C_*^\Delta(M\times M))$ and define
$\mu_2$ to be $\Delta_!\circ\bar{\varepsilon}$.

With these changes, \fullref{rem8.1} becomes correct (it wasn't
before); see \cite{F} for the proof.

In \fullref{sec9}, \fullref{def10.1} should be restated:  define $G_k$ to be the
subcomplex of $(\Sigma^{-m}C_*M)^{\otimes k}$ consisting of elements $x$
for which both $\Sigma^{mk}\bar{\varepsilon}_k(x)$ and 
$\Sigma^{mk}\bar{\varepsilon}_k(\partial x)$ are in general position with
respect to all generalized diagonal maps.

The diagram in \fullref{9.2}(b) should be replaced by
\[
\xymatrix{
\Sigma^{-m_1}C^f_*M_1\otimes 
\Sigma^{-m_2}C^f_*M_2
\ar[r]^{\bar{\varepsilon}}
\ar[d]_{f_!\otimes g_!}
&
\Sigma^{-(m_1+m+2)}C_*^{f\times g}(M_1\times M_2)
\ar[d]^{(f\times g)_!}
\\
\Sigma^{-n_1}C^f_*N_1\otimes 
\Sigma^{-n_2}C^f_*N_2
\ar[r]^{\bar{\varepsilon}}
&
\Sigma^{-(n_1+n+2)}C_*^{f\times g}(N_1\times N_2)
}
\]
 
At the end of \fullref{sec9}, the definition of $G_R$ should be
\[
\bar{\varepsilon}^{-1}_{k'}\circ (R^*)_!\circ 
\bar{\varepsilon}_k.
\]

In \fullref{sec10}, \fullref{10.05} should say that the inclusion
\[
G_{k+l}\hookrightarrow
(\Sigma^{-m}C_*M)^{k+l}
\cong
(\Sigma^{-m}C_*M)^k
\otimes
(\Sigma^{-m}C_*M)^l
\]
has its image in $G_k\otimes G_l$. Now define 
\[
\xi_{k,l}:G_{k+l}\to G_k\otimes G_l
\]
to be the inclusion provided by \fullref{10.05}.

\end{document}